%% file: preprint_arxiv.tex
\documentclass[a4paper, 10pt, DIV10, headinclude=false, footinclude=false]{scrarticle}

\usepackage[utf8]{inputenc}
\usepackage[T1]{fontenc}
\usepackage[english]{babel}

\usepackage[colorlinks=true,linkcolor=blue,citecolor = black ]{hyperref}
\usepackage[pdftex]{graphicx}
\usepackage{subcaption}
\usepackage{latexsym}
\usepackage{amsmath,amssymb,amsthm}
\usepackage{esint} 
\usepackage{aligned-overset}
\usepackage{multirow}
\usepackage{bigdelim}
\usepackage{empheq} 
\usepackage{enumitem}
\usepackage{tikz}
\usepackage{pgfplots}
\usepackage{bm}
\usepackage{xfrac}
\usepackage{wasysym}

\pgfplotsset{compat=newest} 
\usepackage{pgfplotstable}
\usetikzlibrary{positioning,shapes.multipart}
\usepgfplotslibrary{groupplots,dateplot}
\usetikzlibrary{patterns,shapes.arrows}

\newtheorem{Satz}{Theorem}[section]
\newtheorem{Lemma}[Satz]{Lemma}	
\newtheorem{Corollary}[Satz]{Corollary}	
 
\newtheorem{Annahme}[Satz]{Assumption}  

\theoremstyle{definition}
\newtheorem{Definition}[Satz]{Definition} 

\theoremstyle{remark}	 
\newtheorem{Bsp}[Satz]{Example}	  

\newtheorem{Bemerkung}[Satz]{Remark}                  
\numberwithin{equation}{section}

\newcommand{\bspeqqed}{\tag*{\llap{$\Diamond$}}}
\newcommand{\assqed}{\hfill$\ocircle$}
\newcommand{\asseqqed}{\tag*{\llap{$\ocircle$}}}

\newlength\figureheight
\newlength\figurewidth

\newcommand{\R}{\mathbb{R}} 
\newcommand{\N}{\mathbb{N}} 

\newcommand{\timeRscheme}{IMEX }

\newcommand{\iMP}{implicit midpoint }

\newcommand{\cell}[1]{Y_{#1}}

\newcommand{\cellKj}[1]{Y_{#1,K,j}}

\newcommand{\identity}{\operatorname{I}}
\newcommand{\sgn}{\operatorname{sgn}}

\newcommand{\ritzH}{\mathcal{P}_H^{\hilSpac}}

\newcommand{\ritzV}{\mathcal{P}_H^{ \solSpac}}
\newcommand{\ritzblock}{\mathcal{P}_H}
\newcommand{\assOp}{\mathcal{A}}
\newcommand{\assdisOp}{\assOp_H}
\newcommand{\assdampOp}{\mathcal{B}}
\newcommand{\assdisdampOp}{\assdampOp_H}
\newcommand{\interpol}{\mathcal{I}_H}
\newcommand{\interpolblock}{\mathcal{J}_H}
\newcommand{\diffblock}{\mathcal{S}}
\newcommand{\diffdisblock}{S_H}
\newcommand{\diffblockrestrict}{S}

\newcommand{\resOpdis}{\widehat{\mathcal{R}}}
\newcommand{\resplusOpdis}{\resOpdis_+}

\newcommand{\resminusOpdis}{\resOpdis_-}
\newcommand{\resminusOpinvdis}{\resOpdis_-^{-1}}
\newcommand{\resplusminusOpdis}{\resOpdis_\pm}

\newcommand{\nonlin}{\mathcal{G}}

\newcommand{\nonlindis}{\nonlin_H}

\newcommand{\nonlinblock}{F}
\newcommand{\nonlinblockH}{\nonlinblock_H}
\newcommand{\nonlineps}{\nonlin^{\veps}}

\newcommand{\nonlinblockntilde}{\widetilde{\nonlinblock}^{n}}
\newcommand{\nonlinblocktilde}[1]{\widetilde{\nonlinblock}^{#1}}
\newcommand{\nonlinblockminustilde}[1]{\widetilde{\nonlinblock}^{#1-1}}
\newcommand{\nonlinblockhalbetilde}[1]{\widetilde{\nonlinblock}^{#1+\frac{1}{2}}}
\newcommand{\nonlinblockhalbeminustilde}[1]{\widetilde{\nonlinblock}^{#1-\frac{1}{2}}}

\newcommand{\nonlinblockdis}{\nonlinblock_H}
\newcommand{\nonlinblocknhalbetildeH}{\widetilde{\nonlinblock}_H^{n+\frac{1}{2}}}

\newcommand{\nonlinblockntildeH}{\widetilde{\nonlinblock}_H^{n}}

\newcommand{\solSpac}{\mathbf{V}}

\newcommand{\soldualSpac}{\mathbf{V}^{\star}}
\newcommand{\disSpace}{\mathbf{W}_H}
\newcommand{\soldisSpac}{\mathbf{V}_H}
\newcommand{\hilSpac}{\mathbf{H}}
\newcommand{\solSpacblock}{\mathbf{Y}}
\newcommand{\hildualSpac}{\mathbf{H}^{\star}}
\newcommand{\hildisSpac}{\mathbf{H}_H}
\newcommand{\hilSpacblock}{\mathbf{X}}
\newcommand{\soldisSpacblock}{\mathbf{X}_H}
\newcommand{\refSpac}{\mathbf{Z}^\solSpac}
\newcommand{\refSpacblock}{\mathbf{Z}}
\newcommand{\domainblock}{D(\diffblockrestrict)}
\newcommand{\Linf}[1]{L^{\infty}([0,T];#1)}
\renewcommand{\L}[2]{L^{#1}([0,T];#2)}

\newcommand{\Lspace}[1]{L^{#1}(\Omega)}
\newcommand{\Lohnespace}[1]{L^{#1}}

\newcommand{\Hzerospace}[1]{H_0^{#1}(\Omega)}

\newcommand{\Hspace}[1]{H^{#1}(\Omega)}
\newcommand{\Hohnespace}[1]{H^{#1}}

\newcommand{\Vmicro}{V_h^q}

\newcommand{\energynorm}{d}
\newcommand{\energydisnorm}{d_H}
\newcommand{\energytnorm}{\solSpac}
\newcommand{\energytdisnorm}{\soldisSpac}
\newcommand{\energynormdiff}{\Delta\solSpac}
\newcommand{\hilnorm}{\hilSpac}
\newcommand{\hilnormdiff}{\Delta\hilSpac}
\newcommand{\hildisnorm}{\hildisSpac}

\newcommand{\hilskalarprodukt}[2]{(#1,#2)_\hilnorm}
\newcommand{\hilskalarproduktdiff}[2]{(#1,#2)_{\Delta\hilnorm}}
\newcommand{\hilbigskalarprodukt}[2]{\bigl(#1,#2\bigr)_\hilnorm}

\newcommand{\hildisskalarprodukt}[2]{(#1,#2)_{\hildisnorm}}

\newcommand{\dualbracket}[2]{\langle#1,#2\rangle_{\soldualSpac \times \solSpac}}

\newcommand{\norm}[1]{\|#1\|}
\newcommand{\normX}[1]{\norm{#1}_\hilSpacblock}
\newcommand{\normXH}[1]{\norm{#1}_{\soldisSpacblock}}

\newcommand{\normOperatorXH}[1]{\norm{#1}_{\soldisSpacblock \leftarrow \soldisSpacblock}}

\newcommand{\solvec}{x}

\newcommand{\solveczero}{\solvec^{0}}
\newcommand{\solvecn}{\solvec^{\,n}}

\newcommand{\solvecitilde}[1]{\tilde{\solvec}^{#1}}

\newcommand{\solvecnp}{\solvec^{\,n+1}}

\newcommand{\solvechalbn}{\solvec^{\,n+\frac{1}{2}}}

\newcommand{\solvechalbtilde}[1]{\tilde{\solvec}^{\,#1+\frac{1}{2}}}

\newcommand{\soldisvec}{\solvec_H}
\newcommand{\soldisveczero}{\solvec_H^{0}}
\newcommand{\solvectildedis}{\tilde{\solvec}_H}

\newcommand{\solvecndis}{\soldisvec^{\,n}}
\newcommand{\solvecntildedis}{\solvectildedis^n}

\newcommand{\solvechaldis}[1]{\soldisvec^{\,#1+\frac{1}{2}}}

\newcommand{\solvecnH}{\soldisvec^{\,n}}
\newcommand{\solvecntildeH}{\solvectildedis^n}
\newcommand{\solvecitildeH}[1]{\solvectildedis^{#1}}

\newcommand{\solvecnpH}{\soldisvec^{\,n+1}}

\newcommand{\solvecnptildeH}{\solvectildedis^{\,n+1}}
\newcommand{\solvechalbnH}{\soldisvec^{\,n+\frac{1}{2}}}

\newcommand{\solvechalbntildeH}{\solvectildedis^{\,n+\frac{1}{2}}}

\newcommand{\phihj}{\phi_{h,K,j}}
\newcommand{\psihj}{\psi_{h,K,j}}

\newcommand{\phihbarj}{\overline{\phi}_{K_j}}

\newcommand{\nohom}[1]{#1}
\newcommand{\nohomn}[1]{#1^{n}}
\newcommand{\nohomnp}[1]{#1^{n+1}}
\newcommand{\nohomhalbe}[1]{#1^{n + \frac{1}{2}}}
\newcommand{\nohomnull}[1]{#1^{0}}
\newcommand{\nohomdis}[1]{#1_H}
\newcommand{\nohomdisn}[1]{#1^{n}_H}
\newcommand{\nohomdisnull}[1]{#1_H^0}
\renewcommand{\hom}[1]{#1^{hom}}
\newcommand{\homn}[1]{#1^{hom,n}}
\newcommand{\homnull}[1]{#1^{0}}
\newcommand{\homdis}[1]{#1^{hom}_H}
\newcommand{\homdisn}[1]{#1^{hom,n}_H}
\newcommand{\homdisnull}[1]{#1^{0}_H}

\newcommand{\homdisCoeKj}{a^{hom}_{K,j,h}}

\newcommand{\homCoeKj}{a^{hom}_{K,j}}

\newcommand{\hilfAj}{\overline{a}_{K,j}}
\newcommand{\bilinear}[1]{d_{#1}}
\newcommand{\bilinearhom}{d^{hom}}
\newcommand{\bilinearbeta}[1]{b_{#1}}

\newcommand{\Coeffeps}{a^{\veps}}

\newcommand{\abstractspaceError}[1]{E_{H,#1}}
\newcommand{\errorH}{e^n_{\interpolblock}}
\newcommand{\remainderdiffblock}{R_H}
\newcommand{\remainderlipschitz}{r_H}
\newcommand{\errorfunc}{\mathbf{E}_\tau}

\newcommand{\errornonlinear}{\Delta \nonlinblock}

\newcommand{\errornH}{e_H^n}
\newcommand{\errornpH}{e_H^{n+1}}
\newcommand{\errornhalbH}{e_H^{n + \frac{1}{2}}}
\newcommand{\erroriH}{e_H^i}

\newcommand{\errorihalbH}{e_H^{i+ \frac{1}{2}}}

\newcommand{\errornonlinearnH}{\errornonlinear_H^{n}}
\newcommand{\errornonlinearnhalbeH}{\errornonlinear_H^{n + \frac{1}{2}}}
\newcommand{\errornonlineariH}{\errornonlinear_H^{i}}
\newcommand{\errornonlinearihalbeH}{\errornonlinear_H^{i + \frac{1}{2}}}
\newcommand{\abstractfullError}[1]{E_{H,\tau,#1}}

\newcommand{\constemb}{C_{\hilSpac,\solSpac}}
\newcommand{\constembdis}{C_{\hildisSpac,\soldisSpac}}
\newcommand{\constsp}{c_{\energynorm}}
\newcommand{\constspdis}{\widehat{c}_{\energynorm_H}}
\newcommand{\constqm}{c_{qm}}
\newcommand{\constbeta}{c_{\beta}}
\newcommand{\constBeta}{C_{\beta}}

\newcommand{\constqmS}{c^{S}_{qm}}
\newcommand{\constqmdis}{\widehat{c_{qm}}}
\newcommand{\constqmdisS}{c^{S_H}_{qm}}
\newcommand{\constdisH}{C_{\hilSpac}}
\newcommand{\constdisV}{C_{\solSpac}}
\newcommand{\constdisX}{C_{\hilSpacblock}}
\newcommand{\constlipschitz}{L_{\rho}}
\newcommand{\constlipschitzdis}{L_{\rho_H}}
\newcommand{\constcoercive}{\lambda}

\newcommand{\constbounded}{\Lambda}

\newcommand{\constpoincare}{c_{\Omega}}
\newcommand{\constinterpol}{C_{\refSpac}}
\newcommand{\constinterpolblock}{C_{\interpolblock}}
\newcommand{\constremainder}{C_{R}}
\newcommand{\constdefects}{C_{H}}
\newcommand{\constderivative}[1]{M_{#1}}
\newcommand{\constSderivate}{M_{\diffblockrestrict}}

\newcommand{\massmatrix}{\mathbf{M}_H}
\newcommand{\stiffmatrix}{\mathbf{A}_H}
\newcommand{\dampmatrix}{\mathbf{B}_H}
\newcommand{\nonlinearvec}{\mathbf{G}_H}
\newcommand{\loadvec}{\mathbf{f}_H}
\newcommand{\solvecfem}{\bm{\mu}_H}
\newcommand{\solprimevecfem}{\bm{\nu}_H}

\newcommand{\defectCN}[1]{\delta_{\text{MP}}^{#1}}

\newcommand{\defectstagesCN}[1]{D_{\text{MP}}^{#1}}
\newcommand{\defectstagesCNdiff}[1]{\Delta D_{\text{MP}}^{#1}}

\newcommand{\defectIMEXRH}[1]{\delta_{\text{IMEX},H}^{#1}}
\newcommand{\defectIMEXRtotalH}[1]{\delta_{\text{IMEX},H}^{#1}}
\newcommand{\defectstagesIMEXRH}[1]{D_{\text{IMEX},H}^{#1}}

\newcommand{\defectstagesIMEXRdiffH}[1]{\Delta D_{\text{IMEX},H}^{#1}}
\newcommand{\defectspaceMP}[1]{\delta_H^{#1}}
\newcommand{\defectspaceMPstage}[1]{d_H^{#1}}

\newcommand{\e}{\operatorname{e}}
\renewcommand{\d}[1]{\, \operatorname{d}\!#1}

\newcommand{\veps}{\varepsilon}
\newcommand{\coloneqq}{:=}
\newcommand{\tauhalb}{\frac{\tau}{2}}
\newcommand{\sumi}[2]{\sum^{#2}_{i=#1}}
\newcommand{\suminullnminuseins}{\sum^{n-1}_{i=0}}

\usepackage{lipsum}

\newenvironment{abstr}[1]{ \vspace{.05in}\footnotesize
	\parindent .2in
	{\upshape\bfseries #1. }\ignorespaces}{\par\vspace{.1in}}
	\newenvironment{Abstract}{\begin{abstr}{Abstract}}{\end{abstr}}
	\newenvironment{keywords}{\begin{abstr}{Key words}}{\end{abstr}}
	\newenvironment{AMS}{\begin{abstr}{AMS subject classifications}}{\end{abstr}}
	\makeatletter
\newcommand{\mylabel}[2]{#2\def\@currentlabel{#2}\label{#1}}
\makeatother

\allowdisplaybreaks[4]

\begin{document}  

	\title{Error analysis of an implicit-explicit time discretization scheme for semilinear wave equations with application to multiscale problems%
	\thanks{Funded by the Deutsche Forschungsgemeinschaft (DFG, German Research Foundation) – Project-ID 258734477 – SFB 1173. BV also acknowledges funding by the Deutsche Forschungsgemeinschaft (DFG, German Research Foundation) under project number 496556642 and under Germany's Excellence Strategy –  EXC-2047/1 – 390685813.}
	}
	\author{Daniel Eckhardt\footnotemark[2] \and Marlis Hochbruck\footnotemark[2] \and Barbara Verf\"urth\footnotemark[3]}
	\date{}
	\maketitle
	
	\renewcommand{\thefootnote}{\fnsymbol{footnote}}
	\footnotetext[2]{Institut für Angewandte und Numerische Mathematik, Karlsruher Institut für Technologie, Englerstr. 2, D-76131 Karlsruhe}
	\footnotetext[3]{Institut für Numerische Simulation, Universit\"at Bonn, Friedrich-Hirzebruch-Allee 7, D-53115 Bonn}
	\renewcommand{\thefootnote}{\arabic{footnote}}

	\begin{Abstract} 
		We present an implicit-explicit (IMEX) scheme for  semilinear  wave  equations with strong damping. By treating the  nonlinear, nonstiff term explicitly and the linear, stiff part implicitly, we obtain a method which is not only unconditionally stable but also highly efficient. Our main results are error bounds of the full discretization in space and time for the IMEX scheme combined with a general abstract space discretization. As an application, we consider the heterogeneous multiscale method for wave equations with highly oscillating coefficients in space for which we show spatial and temporal convergence rates by using the abstract result. 
	\end{Abstract}



\section{Introduction}
For numerically solving  partial differential equations (PDEs) in time, a variety of  different methods  have been investigated in the literature, broadly categorized as implicit or explicit schemes. Explicit schemes are popular since they are simple to implement, computationally efficient, but they suffer from step-size restrictions, so-called CFL conditions. On the other hand, implicit schemes involve solving a system of nonlinear equations at each time step, resulting in improved stability and allowing for larger time steps compared to explicit methods. However, they require a higher computational effort. To combine the advantages of both strategies, we construct and numerically analyze an implicit-explicit (IMEX) method in which the CFL condition only depends  on  non-stiff terms. Besides many applications, e.g., air pollution models \cite{Verwer1996Blom}, hydrodynamics and atmospheric dynamics \cite{Gardner2018atmospheric,Kadioglu2010hydro}, fluid-structure interaction \cite{Pedro2015Brandastrucdynamics}, there is a well-developed theory of IMEX-Runge Kutta (IMEX-RK) \cite{Ascher1997Ruuth,Boscarino2007,Koto2008,Layton2015Li,Li2021Quan} and IMEX multistep methods \cite{Akrivis1999Crouzeix,Ascher1995Ruuth,Huang2022Shen,Hundsdorfer2007Ruuth} for ordinary differential equations and parabolic problems. 
For the wave equation IMEX-RK \cite{Arun2021Gupta}, Crank-Nicolson/Leapfrog (CNLF) schemes \cite{Leibold2021Hochbruck,Layton2012Trenchea} (note that these are not equivalent) and a $\theta$-scheme/Leapfrog combination \cite{Methenni2023Imperiale} are considered. 
As in \cite{Leibold2021Hochbruck}, we consider numerical solutions of semilinear wave equations in a general Hilbert space $\hilSpac$ of the following form 
\begin{align}
    \nohom{u}^{\prime\prime }+ \assOp u & + \assdampOp \nohom{u}^{\prime} =  \nonlin(u,u^{\prime}) + f. \label{intro::EE}
\end{align}
Here, $\assOp $ and $  \assdampOp$ are linear, possibly unbounded operators, representing the wave propagation and (strong) linear damping of the system, respectively. The nonlinear operator $\nonlin$ and the function $f$ denote additional (nonlinear) source and damping terms. In \cite{Leibold2021Hochbruck} a rigorous error analysis for the CNLF applied to \eqref{intro::EE} was presented  if $\nonlin(u,u^\prime) = \nonlin(u)$. 
However, if $\nonlin$ depends on $u^\prime$, this scheme is not very attractive, as it leads to a system of nonlinear equations. In addition, \cite{Leibold2021Hochbruck} requires boundedness conditions on $\assdampOp$ which are not be fulfilled in our application. We therefore propose a different method based on the implicit and explicit midpoint rule applied to the first-order system of \eqref{intro::EE}. Our scheme leads to a system of \textit{linear} equations also in the general case and can be formulated for a larger class of operators $\assdampOp$.  Even with these greater application possibilities, the computational effort of the implicit/explicit midpoint rule is about the same as for the CNLF in the case $\nonlin(u,u^\prime) = \nonlin(u)$ (see the discussion in Section \ref{num_example}). Although the scheme has already been considered in the literature \cite{Ascher1995Ruuth}, to the best of our knowledge there is no error analysis for (nonlinear) wave equations.  We provide a proof of second-order convergence, which requires significantly different arguments than in \cite{Leibold2021Hochbruck}. A central aspect of the proof is to consider the first-order formulation and to interpret the scheme as a perturbation of the implicit midpoint rule. This view allows us to bound the error terms of the IMEX method by terms of the implicit midpoint rule. Building on the unified error analysis from \cite{Hipp2019}, our main result is a full discretization bound for the IMEX scheme combined with an abstract space discretization method.

As an application we consider the semilinear wave equation with damping
\begin{equation}
    \begin{aligned}
        \partial_{tt} u^{\veps} - \nabla\cdot (\Coeffeps(x)\nabla u^{\veps}) - \nabla\cdot ( \beta(x) \nabla \partial_tu^{\veps})    =   \nonlin(x,\partial_{t}u^{\veps})     + f  \label{intro::application}
    \end{aligned}
\end{equation}
in a bounded Lipschitz domain $\Omega$ and with Dirichlet boundary conditions and  initial values. The parameter $\veps$ describes small scale effects, e.g., fast oscillations, which appear in the coefficient $\Coeffeps$ and therefore also in the solution $u^\veps$. This in connection with Finite Elements (FE)  can lead to a high computational effort, since a grid width in the order of magnitude of $\veps$ must be chosen.  On the other hand, in applications we are often  only interested in the macroscopic behavior of $u^\veps$. We therefore propose a multiscale method based on the Heterogeneous Multiscale Method (HMM) introduced in \cite{Abdulle2012Engquist}. 
The idea of the HMM is to derive a model related to \eqref{intro::application}, which is independent of $\veps$. Missing microscopic data is collected by solving microscopic problems on small subdomains and only where the data is needed.
Here we make use of a convergence result for $\veps \to 0$ (see \cite{Nguetseng2010}) which leads to a homogenized version of \eqref{intro::application} with an $\veps$-independent coefficient $\hom{a}$. This coefficient can be calculated explicitly if scale separation and $\veps$-periodicity of $a^\veps$ is assumed, by solving  so-called cell problems.
  For the space discretization we consider FE for the homogenized problem of \eqref{intro::application}, where we can now choose a much coarser mesh than if we would apply FE directly for \eqref{intro::application}, since the homogenized problem does not contain any microscopic features. We show that the HMM fits into the abstract framework and thus the error estimates can be applied, which extends the known results for the HMM discretization from linear \cite{Abdulle2011Grote,Abdulle2014} to semilinear wave equations.  

The paper is organized as follows. In Section \ref{generalSet} we introduce the abstract setting, present a wellposedness result, the time integration scheme, and the error result for the semidiscrete equation. The full discretization error is investigated in Section \ref{fulldis} where we combine the IMEX scheme with an abstract space discretization. Finally, we derive the fully discrete finite element heterogeneous multiscale methode for semilinear wave equations with damping in Section \ref{Application} and prove convergence rates using the abstract results. These theoretical findings are demonstrated by
numerical results presented in Section \ref{num_example}. 
\section{Time Integration \label{generalSet}}
In this section, we consider time-integration methods for the first-order model equation \eqref{intro::EE}.  We start with the general framework and the wellposedness result. Secondly, we construct an IMEX scheme  motivated by the observation that the operator $\nonlin$ is often non-stiff in practice. 

\subsection{Wellposedness of the continuous problem}\label{subsec:wellposed}

    {Let $\hilSpac$ be a Hilbert space with inner products $(\cdot,\cdot)_{\hilSpac}$ and $\solSpac \subset \hilSpac$ be a  densely embedded subspace. The corresponding norm is denoted by $\|\cdot\|_{\hilnorm}$.
        We consider the variational problem
        \begin{equation}
            \begin{aligned}
                \label{Setting::weak}
                \hilbigskalarprodukt{u^{\prime \prime}}{w} + \bilinear{}(\nohom{u}, w) &+  \bilinearbeta{}(  \nohom{u}^{\prime},{w} ) =  \hilskalarprodukt{\nonlin(u, u^\prime)}{w} + \hilbigskalarprodukt{f(t)}{w} \\ 
                u(0) &= \nohomnull{u} , \qquad   u^\prime(0) = \nohomnull{v}
            \end{aligned}
        \end{equation}
        for all $w \in \solSpac$, $t\geq 0$,
        and impose the following conditions on $\bilinear{}, \bilinearbeta{}$, $\nonlin$ and $f$ :
        \begin{Annahme}[Forms, nonlinearity, source term]  
            \label{Setting::assumption} \phantom{42}
            \begin{enumerate}[label=\textnormal{(\arabic*)}]
                \item[(i) ] The bilinear form $\energynorm: \solSpac \times  \solSpac \to \R$ is symmetric and continuous. In addition, there exist constants $\constsp$ and  $\constemb$ such that
                \begin{align*}
                    (\cdot,\cdot)_\energytnorm = \energynorm(\cdot,\cdot) + \constsp \, (\cdot,\cdot)_\hilnorm
                \end{align*}
                is an inner product on $\solSpac$ and the embedding $\solSpac \hookrightarrow \hilSpac$ is continuous, i.e., 
                \[
                \|v\|_{\hilnorm} \leq \constemb\|v\|_{\energytnorm} \qquad \text{ for all } v \in \solSpac.
                \]
                \item[(ii) ] The bilinear form $ \bilinearbeta{}: \solSpac \times \solSpac \to \R$  is continuous and quasi-monotone, i.e., there exists a constant $\constqm \geq 0$ such that:
                \begin{align*}
                    \bilinearbeta{}({ v },{v}) \geq -\constqm \|v\|^2_{\hilnorm}.
                \end{align*}
                \item[(iii) ] The nonlinearity $\nonlin : \solSpac \times \hilSpac  \to \hilSpac$ is locally Lipschitz continuous , i.e, for all $\rho > 0$ there exists a Lipschitz constant $L_{\rho}$ such that for all $v_1,v_2 \in \solSpac$ and $w_1,w_2 \in \hilSpac$ with  $ \max \{\|[v_1,w_1]\|_{\energytnorm\times\hilnorm},\|[v_2,w_2]\|_{\energytnorm\times\hilnorm} \} \leq \rho$ it holds
                \begin{align*}
                    \| \nonlin(v_1,w_1) - \nonlin(v_2,w_2)\|_{\hilnorm} \leq L_{\rho} ( \| v_1 - v_2 \|_{\energytnorm} + \| w_1 - w_2 \|_{\hilnorm}),
                \end{align*}
                \item[(iv) ]  The right-hand side $ f: [0,\infty) \to \hilSpac $ satisfies
                \begin{align*}
                    f \in W_{loc} ^{1,1}([0,\infty);\hilSpac).         
                    \asseqqed
                \end{align*}
            \end{enumerate}
        \end{Annahme}
        
         By (i), $(\cdot,\cdot)_\energytnorm$, induces a norm $\|\cdot\|_{\energytnorm}$. The product space $\solSpac \times \hilSpac$ will be equipped with the norm  $\|\cdot\|_{\energytnorm\times\hilnorm} \coloneqq \|\cdot\|_{\energytnorm} + \|\cdot\|_{\hilnorm}$.  Moreover, we identify $\hilSpac$ with its dual space $\hildualSpac$, resulting in  the following Gelfand triple:
\begin{align*}
    \solSpac \hookrightarrow \hilSpac = \hildualSpac \hookrightarrow \soldualSpac.
\end{align*}

We define the operators $\assOp, \assdampOp \in \mathcal{L}(\solSpac,\soldualSpac)$ associated to $\bilinear{}, \bilinearbeta{}$ :
\begin{align*}
    \dualbracket{\assOp v}{w} &:=  \bilinear{}(v, w), \\
    \dualbracket{\assdampOp v}{w} &:= \bilinearbeta{}(v, w).
\end{align*}
Problem \eqref{Setting::weak} can thus be written as the following evolution equation 
\begin{align}
    \label{Setting::EE}
    \nohom{u}^{\prime\prime }+ \assOp u  + \assdampOp \nohom{u}^{\prime} &=  \nonlin(u,u^{\prime}) + f \\
    u(0) = \nohomnull{u} \quad        & ,  u^\prime(0) = \nohomnull{v}. \notag
\end{align}
\subsection{First-order formulation}
In order to derive a first-order formulation of \eqref{Setting::EE}, we set:  $\nohom{u}^{\prime} = \nohom{v}$,
\begin{equation*}
    \begin{aligned}
        \solvec =\left[ \begin{array}{rrr}
            \nohom{u} \\
            \nohom{v} \\
        \end{array}\right], \quad
        \diffblock\solvec =\left[ \begin{array}{rrr}
            \nohom{v}       &                                         \\
            -\assOp \nohom{u} & \hspace*{-0.3cm}-  \assdampOp \nohom{v}  \\
        \end{array}\right],
        \quad
        \nonlinblock(\solvec) = \left[
        \begin{array}{rrr}
            0           &         \\
            \nonlin(\nohom{u},\nohom{v})& \hspace*{-0.3cm}+f  \\
        \end{array}\right],
        \quad
        \solveczero =\left[
        \begin{array}{rrr}
            \nohomnull{u} \\
            \nohomnull{v} \\
        \end{array}\right].
    \end{aligned}
\end{equation*}
Then \eqref{Setting::EE} reads as follows: \\ Seek $\solvec  \in \solSpacblock := \solSpac \times \solSpac$, such that
\begin{equation}
    \label{Setting::fof}
    \begin{aligned}
        \solvec^{\prime} - \diffblock\solvec & =  \nonlinblock(\solvec) \\
        \solvec(0)                                 & = \solveczero.
    \end{aligned}
\end{equation}
To obtain an evolution equation in $\hilSpacblock :=\solSpac \times \hilSpac $ 
we consider the following restriction of the  operator $\diffblockrestrict$:
\begin{align*}
    \diffblockrestrict: \domainblock \to \hilSpacblock, \qquad y \to \diffblockrestrict y = \diffblock y \qquad \text{on } \domainblock & =  \{y \in \solSpacblock \mid \diffblock y  \in \hilSpacblock\} \\ &= \{(u,v) \in \solSpac \times \solSpac \mid \assOp u + \assdampOp v\in \hilSpac\}.
\end{align*}
\begin{Lemma}
    \label{Settig::lemma:fo_qmono}
    The operator $\diffblockrestrict: \domainblock \to \hilSpacblock$ is the generator of $C_0$-semigroup, which satisfies \[ \|\e^{t \diffblockrestrict}\|_{\hilSpacblock \leftarrow \hilSpacblock} \leq \e^{\constqmS t}\] with   $\constqmS \coloneqq\frac{1}{2}\constsp \constemb + \constqm$.
\end{Lemma}
\begin{proof}
    $\,$ This follows by a combination of Lemma 4.2 (with $\alpha = 1$), Lemma 2.3 and Theorem 2.4 in \cite{Hipp2019}. 
\end{proof}
The well-posedness result can be found in \cite[Lemma 3.2]{Leibold2020}.
\begin{Corollary}
    \label{Setting::Corollary:wellposedness}
    Let Assumption 	$\ref{Setting::assumption}$ be satisfied and $\nohomnull{u}$, $\nohomnull{v}$  $ \in\solSpac$ with $\assOp \nohomnull{u} + \assdampOp \nohomnull{v} \in \hilSpac$. Then there exist a maximal existence time $t^{\star}(\nohomnull{u},\nohomnull{v}) \in (0, \infty]$ such that for all $T < t^{\star}(\nohomnull{u},\nohomnull{v})$, \eqref{Setting::EE} has a unique solution $\nohom{u} \in W^{2,\infty}([0,T];\hilSpac) \cap C^1([0,T],\solSpac)$ which satisfies $\assOp \nohom{u} + \assdampOp \nohom{u}^{\prime} \in C([0,T];\hilSpac)$.
\end{Corollary}

\subsection{The \timeRscheme scheme}
The \iMP scheme applied \eqref{Setting::fof} can be written in a half step formulation as
\begin{subequations}
    \label{Setting::FoF}
    \begin{alignat}{2}
        \solvechalbn &= \solvecn + \frac{\tau}{2}\diffblockrestrict \solvechalbn  + \frac{\tau}{2}\nonlinblock(\solvechalbn), \label{Setting::FoF:1}\\
        \solvecnp &= \solvecn + \tau \diffblockrestrict \solvechalbn + \tau \nonlinblock(\solvechalbn),\label{Setting::FoF:2}
    \end{alignat}
\end{subequations}
while explicit midpoint rule uses 
\begin{equation}
    \label{IMEX::explicitMP}
    \begin{aligned}
        \solvechalbn &= \solvecn + \frac{\tau}{2}\diffblockrestrict \solvecn  + \frac{\tau}{2}\nonlinblock(\solvecn)
    \end{aligned}
\end{equation}
instead of \eqref{Setting::FoF:1}.
Combining \eqref{Setting::FoF:1} and \eqref{IMEX::explicitMP} leads to the following linearly implicit scheme:
\begin{subequations}
    \label{IMEX::IMEXR}
    \begin{alignat}{2}
        \solvechalbn &= \solvecn + \frac{\tau}{2}\diffblockrestrict \solvechalbn + \frac{\tau}{2}\nonlinblock(\solvecn), \label{IMEX::IMEXR1}\\
        \solvecnp &= \solvecn + \tau \diffblockrestrict \solvechalbn + \tau \nonlinblock(\solvechalbn).\label{IMEX::IMEXR2}
    \end{alignat}
\end{subequations}
Hence the $\diffblockrestrict$ operator is treated implicitly and the nonlinear operator $\nonlinblock$ explicitly. This results in  second-order error bounds for the \timeRscheme scheme. 
\begin{Satz}
    \label{IMEX::Theorem}
    Let $u \in C^4([0,T]; \hilSpac) \cap C^3([0,T];\solSpac)$ be the solution of \eqref{Setting::EE} with  $\assOp \nohom{u} + \assdampOp u^\prime \in C^2([0,T];\hilSpac)$ and $\solvecn = [\nohomn{u}, \nohomn{v}]^{\intercal}$, $t_n= \tau n \in [0,T]$, be the approximations obtained by the \timeRscheme scheme \eqref{IMEX::IMEXR}. For
    \begin{align*}
        \rho  = 2(\|\nohom{u}\|_{\L{\infty}{\solSpac}} + \|\nohom{u}^{\prime}\|_{\L{\infty}{\hilSpac}}),
    \end{align*}
    there exists $\tau^\star>0$ such that for all $\tau < \tau^\star$ and all $t_n \in [0,T]$ 
    \[ \max_{t_n \leq T }\|\nohomn{u}\|_{\energytnorm}  \leq \rho \quad \text{ and } \quad \max_{t_n \leq T }  \|\nohomn{v}\|_{\hilnorm} \leq \rho.\]
    Furthermore, the error bound
    \begin{align*}
        \| \nohomn{u} - \nohom{u}(t_n)\|_{\solSpac} + \| \nohomn{v} - \nohom{u}^{\prime}(t_n) \|_{\hilSpac} \leq C t_n \e^{ \bigl(\constqmS +  \constlipschitz (2+ \tauhalb \constlipschitz )  \bigr) t_n   }  \tau^2
    \end{align*}
    holds true with a constant $C$ which depends on derivates of $u$ and  $L_{\rho}$, $T$ but is independent of $n$ and $\tau$.
\end{Satz}
Since the proof works similarly as the more complicated proof of Theorem \ref{fulldis::TheoremIMEX} of the full discretization, we omit it here,  cf.  Remark \ref{fulldis::remark} for more details.

\section{Full discretization \label{fulldis}}
In this section we analyze the full discretization of \eqref{Setting::EE} using the \timeRscheme scheme \eqref{IMEX::IMEXR} as time stepping method and an abstract space discretization. We first recall some results from \cite{Leibold2020}. 
\subsection{Space discretization}
Let $(\hildisSpac)_H \subset \solSpac$  be a family of finite-dimensional spaces, e.g., a finite element space.
For these subspaces we can formulate a discrete version of \eqref{Setting::EE} as
\begin{equation}
    \begin{aligned}
        \label{Spacedis::weakdis}
        \hildisskalarprodukt{\nohomdis{u}^{\prime \prime}}{\nohomdis{w}} & + \energydisnorm(\nohomdis{u}, \nohomdis{w}) + \bilinearbeta{H}(\nohomdis{u}^\prime, \nohomdis{w}) \\ & +  \hildisskalarprodukt{\nonlindis(\nohomdis{u} ,\nohomdis{u}^\prime)}{\nohomdis{w}} =  \hildisskalarprodukt{ f_H(t)}{\nohomdis{w}}
    \end{aligned}
\end{equation}
for $t \geq 0 $,  $\nohomdis{w} \in \hildisSpac$ with $\nohomdis{u}(0) = \nohomdisnull{u} ,\nohomdis{v}(0) = \nohomdisnull{v}$. Here, $ (\cdot, \cdot)_{\hildisnorm}, \energydisnorm, b_H, \nonlindis,\nohomdisnull{u},\nohomdisnull{v}$ are approximations to $(\cdot, \cdot)_{\hilnorm} , \energynorm, b, \nonlin, \nohomnull{u}, \nohomnull{v}$. 	Similar to the continuous case in Assumption~\ref{Setting::assumption}, we assume the following conditions:
    \begin{Annahme}[Discrete setting]
    \label{Spacedis::assumptions}
    \phantom{-} \hspace*{0.1cm}
    \begin{enumerate}[label=\textnormal{(\arabic*)}]
        \item[(i) ] The bilinear form $\energydisnorm: \hildisSpac \times  \hildisSpac \to \R$ is symmetric and continuous, and 
        \begin{align*}
            (\cdot,\cdot)_{\energytdisnorm} =
            \energydisnorm(\cdot,\cdot) + \constsp (\cdot,\cdot)_{\hildisnorm}
        \end{align*}
        is a scalar product on $\hildisSpac$.   
        The approximation space equipped with the  scalar product  $(\cdot,\cdot)_{\energytdisnorm} $
         will be denoted by $\soldisSpac$ with induced norm $\| \cdot \|_{\energytdisnorm}$.
         In addition, we assume that there is a constant $\constembdis > 0$ such that
         \[ 
         \|v_H\|_{\hildisnorm} \leq \constembdis \|v_H\|_{\energytdisnorm} \quad  \text{ for all } v_H \in \soldisSpac.
         \]
        \item[(ii) ] The bilinear form $ \bilinearbeta{H}: \soldisSpac \times \soldisSpac \to \R$  is continuous and quasi-monotone  w.r.t.~the $\|\cdot\|_{\hildisnorm}$-norm with constant $\constqmdis$.
        \item[(iii) ] The nonlinearity $\nonlindis : \soldisSpac \times \hildisSpac \to \hildisSpac$ is locally Lipschitz continuous with Lipschitz constant $\constlipschitzdis$ and w.r.t.~the $\|\cdot\|_{\energytdisnorm \times \hildisnorm }$-norm. 
        \item[(iv) ] The right-hand side $f_H :[0,\infty) \to \hildisSpac$ satisfies
        \begin{align*}
            f_H \in W_{loc} ^{1,1}([0,\infty);\hildisSpac). \asseqqed
        \end{align*}
    \end{enumerate}
\end{Annahme}

Note that in (i) a different Garding constant than $\constsp$ could be used. Then one has to set $\constsp$ as the maximum of the constants arising in Assumptions~\ref{Setting::assumption} and Assumptions~\ref{Spacedis::assumptions}. 

The discrete associated operators $\assdisOp, \assdisdampOp \in \mathcal{L}(\soldisSpac,\soldisSpac)$ are defined via
\begin{align*}
    \hildisskalarprodukt{ \assdisOp v}{w}  :=\bilinear{H}(v,  w), \\
    \hildisskalarprodukt{ \assdisdampOp v}{w}  :=\bilinearbeta{H}(v,  w).
\end{align*}
This yields the following evolution equation, which is the discrete counterpart of \eqref{Setting::EE}:
\begin{equation}
    \label{Spacedis::EE}
    \begin{aligned}
         \nohomdis{u}^{\prime\prime}+ \assdisOp \nohomdis{u} + \assdisdampOp \nohomdis{u}^{\prime} & = \nonlindis(\nohomdis{u},\nohomdis{u}^{\prime}) + f_H	\qquad \text{ in }	\soldisSpac, \\ 
         u_H(0) = \nohomdisnull{u} \quad        & ,  u_H^\prime(0) = \nohomdisnull{v}.
    \end{aligned}
\end{equation}
To derive error bounds for the abstract spatial discretization, we need interpolation and projection operators for which the following assumptions  hold.
\begin{Annahme}[Operators]
    \label{Spacedis::Operators}
    \phantom{-} \hspace*{0.1cm}
    \begin{enumerate}[label=\textnormal{(\arabic*)}]
        \item[(i) ] \textit{Projection operator:}
        There exist linear operators $\ritzV: \solSpac \to \soldisSpac$ and $\ritzH: \hilSpac \to \hildisSpac$ such that:
        \begin{align*}
            \hildisskalarprodukt{\ritzH v}{w_H} & := \hilskalarprodukt{v}{w_H} \quad   \text{ for all } v \in \hilSpac, w_H \in \hildisSpac \\
            (\ritzV v, w_H)_{\energytdisnorm}      & := (v,  w_H)_{\energytnorm} \quad \;   \text{ for all } v \in \solSpac, w_H \in \soldisSpac
        \end{align*}
        for all $t \in[0,T]$.
        \item[(ii) ] \textit{Interpolation operator:}
        There exists an interpolation operator $\interpol \in \mathcal{L}(\refSpac,\soldisSpac)$ defined on a dense subspace $\refSpac$ of $\solSpac$ with \[\|\interpol \|_{\soldisSpac \leftarrow \refSpac} \leq \constinterpol.  \]
        \item[(iii) ] \textit{Embedding:} There exist constants $\constdisH, \constdisV$ such that
        \begin{align*}
            \|v_H\|_{\hilnorm} \leq \constdisH \|v_H\|_{\hildisnorm}, \qquad \qquad \|v_H\|_{\energytnorm} \leq \constdisV \|v_H\|_{\energytdisnorm}
        \end{align*}
        for all $v_H \in \soldisSpac$.\assqed
    \end{enumerate}
\end{Annahme}
\begin{Bemerkung}
     In \cite{Leibold2020} lift operators are introduced that are not needed here, as we assume that $(\hildisSpac)_H \subset \solSpac$. 
\end{Bemerkung}
In addition,  we have to measure the error between the scalar products $(\cdot,\cdot)_\solSpac$, $(\cdot,\cdot)_\hilSpac$ and the corresponding discrete counterparts.
\begin{Definition}
    For $v_H,w_H \in \soldisSpac$ we define the following errors by the differences of the scalar product:
    \begin{align*}
         \hilskalarproduktdiff{v_H}{w_H} & := \hilskalarprodukt{v_H}{w_H} - \hildisskalarprodukt{v_H}{w_H} \\
         (v_H,w_H)_{\energynormdiff}       & := (v_H,w_H)_{\energytnorm} - (v_H,w_H)_{\energytdisnorm}.            
    \end{align*}
\end{Definition}
Analogously to the continuous case, we now want to write the discrete equation \eqref{Spacedis::EE} in first-order formulation. This allows us to transfer the error result from \cite{Leibold2020}.  \\
For the semidiscrete equation we set: $\nohomdis{u}^{\prime} = \nohomdis{v}$
\begin{align*}
    \soldisvec  =\left[\arraycolsep=0.4pt\def\arraystretch{1.2} \begin{array}{rrr}
        \nohomdis{u} \\
        \nohomdis{v} \\
    \end{array}\right], \; \;
    \diffdisblock \soldisvec =\left[\arraycolsep=0.4pt\def\arraystretch{1.2} \begin{array}{rrr}
        \nohomdis{v}            &                             \\
        -\assdisOp \nohomdis{u} - & \assdisdampOp\nohomdis{v} \\
    \end{array}\right], \; \;
    \nonlinblockdis(\soldisvec) =\left[
        \begin{array}{rrr}
            0           &         \\
            \nonlindis(\nohomdis{u}, \nohomdis{v})& \hspace*{-0.3cm}+f_H  \\
        \end{array}\right], \;                       \soldisveczero =\left[ \arraycolsep=0.4pt\def\arraystretch{1.2}
    \begin{array}{rrr}
        \nohomdisnull{u} \\
        \nohomdisnull{v} \\
    \end{array}\right].
\end{align*}
We can write the discrete equation as  follows \\
Seek $x_H \in \soldisSpacblock \coloneqq \soldisSpac \times \hildisSpac$, such that
\begin{equation}
    \begin{aligned}
        \label{Spacedis::fofdis}
         \soldisvec^\prime - \diffdisblock\soldisvec & =   \nonlinblockdis (\soldisvec)          \\
        \soldisvec(0)                                                                     & = \soldisveczero.
    \end{aligned}
\end{equation}
Local well-posedness for the first-order formuation \eqref{Spacedis::fofdis} and therefore local well-posedness for \eqref{Spacedis::EE} follows from Picard-Lindelöf theorem with maximal existence time of the solution  $t^\star_H(\solveczero_H)$. \\
Then, using the operators from  Assumption \ref{Spacedis::Operators}, the operators for the first-order formulation \eqref{Spacedis::fofdis} are given by the following definition.
\begin{Definition}
    \label{Spacedis::fo_operators}
    We set $\refSpacblock = \solSpac \times \refSpac$ and define the first-order reference operator $\interpolblock: \refSpacblock \to \soldisSpacblock $ and first-order projection operator $\ritzblock: \hilSpacblock \to \soldisSpacblock $ by
    \begin{align*}
        \interpolblock\left[
        \begin{array}{rrr}
            v \\
            w \\
        \end{array}\right]
        : = \left[\arraycolsep=0.4pt\def\arraystretch{1.2}
        \begin{array}{rrr}
            \ritzV v    \\
            \interpol w \\
        \end{array}\right], \qquad
        \ritzblock\left[
        \begin{array}{rrr}
            \nohom{v} \\
            \nohom{w} \\
        \end{array}\right]
        : = \left[\arraycolsep=0.4pt\def\arraystretch{1.2}
        \begin{array}{rrr}
            \ritzV \nohom{v} \\
            \ritzH \nohom{w} \\
        \end{array}\right].
    \end{align*}
\end{Definition}
By Assumption \ref{Spacedis::Operators} we have 
\begin{equation*}
    \begin{aligned}
        \|\interpolblock\|_{\soldisSpacblock \leftarrow \refSpacblock} \leq \constinterpolblock \coloneqq \constdisV + \constinterpol  , \qquad \|\ritzblock\|_{\soldisSpacblock \leftarrow \hilSpacblock} \leq \constdisX \coloneqq \constdisV + \constdisH . 
    \end{aligned}
\end{equation*}

 The reference operator $\interpolblock$ is needed to obtain optimal error bounds w.r.t $H$. 
\subsection{Error bounds for the full discretization}
For the fully discrete scheme, we apply the IMEX scheme to \eqref{Spacedis::fofdis}. This yields
\begin{subequations}
    \label{errorbounds::IMEXR}
    \begin{alignat}{2}
        \solvechaldis{n} &= \solvecnH + \frac{\tau}{2}\diffdisblock \solvechalbnH + \frac{\tau}{2}\nonlinblockH(\solvecnH), 		\label{errorbounds::IMEXR1}\\
        \solvecnpH &= \solvecnH + \tau \diffdisblock \solvechaldis{n} + \tau \nonlinblockH(\solvechaldis{n}). 		\label{errorbounds::IMEXR2}
    \end{alignat}
\end{subequations}

We now present the main result, i.e., error bounds for the fully discrete \timeRscheme scheme. 
\begin{Satz} 
    \label{fulldis::TheoremIMEX}
    Let $u \in C^4([0,T]; \hilSpac) \cap C^3([0,T];\solSpac)$ be the solution of \eqref{Setting::EE} with $\nohom{u},\nohom{u}^{\prime},\nohom{u}^{\prime\prime} \in L^{\infty}([0,T];\refSpac)$ and $\assOp \nohom{u} + \assdampOp u^\prime \in C^2([0,T];\hilSpac)$.  Further, let  $\solvecndis = [\nohomdisn{u} \, \nohomdisn{v}]^{\intercal}$ be the approximations obtained by the fully discrete \timeRscheme scheme \eqref{IMEX::IMEXR} at time $t_n \in [0,T]$. We set
    \begin{align*}
        \rho_H: = \max \bigl\{ (\constdisV + \constinterpol) \|[\nohom{u},{u}^{\prime}]\|_{\L{\infty}{\solSpac \times \refSpac}},\max_{t_n \leq T} \|\nohomn{x}_H\|_{\energytdisnorm \times\hildisnorm}  \bigr\}
    \end{align*}
    and assume that $\tau$ satisfies the step size restriction
    \begin{align}
        \label{fulldis::timestepsizerestriction}
        \tau < \frac{1}{\constqmdisS}.
    \end{align}
    Then, the error bound
    \begin{align}
        \label{fulldis::IMEX:errorboundtotal}
              \| \nohomdisn{u} - \nohom{u}(t_n)\|_{\solSpac} + \| \nohomdisn{v} - \nohom{u}^{\prime}(t_n) \|_{\hilSpac} \leq C t_n \e^{ M t_n   }  (\tau^2 + E_H)
    \end{align}
    holds true with $ M \coloneqq \constqmdisS +  \constlipschitzdis (2+ \tauhalb \constlipschitzdis ) $ and a constant $C$ which depends on derivatives of $u$, but is independent of $n$ and $\tau$.
    Here,
            \begin{subequations} \label{errorbounds::abstractspaceError}
            \begin{equation} \label{eq:EHdef}
                  E_H = E_H(u) =  E_{H,1}+ E_{H,2} + E_{H,3} + E_{H,4}
            \end{equation}
            contains the abstract space discretization error terms $E_{H,i} = E_{H,i}(u,f)$ given by
    \begin{equation}
        \label{errorbounds::abstractspaceError-1-4}
        \begin{aligned}
            E_{H,1} & = \|\nohomdisnull{u} - \ritzV \nohomnull{u}\|_{\energytdisnorm} + \| \nohomdisnull{v} - \interpol\nohomnull{v}\|_{\hildisnorm} + \|\ritzH f - f_{H}\|_{\L{\infty}{\hildisSpac}}                                                               \\
            E_{H,2} & = \max_{\|\phi_{H}\|_{\energytdisnorm} = 1} \| (\interpol \nohom{u}, \phi_H)_{\energynormdiff}\|_{L^{\infty}[0,T]}
            + \max_{\|\psi_{H}\|_{\hildisnorm} = 1} \| (\interpol \nohom{u}, \psi_H)_{\hilnormdiff} \|_{L^{\infty}[0,T]}   \\ &  \qquad +\max_{\|\phi_{H}\|_{\energytdisnorm} = 1} \| (\interpol \nohom{u}^\prime, \phi_H)_{\energynormdiff}\|_{L^{\infty}[0,T]} + \max_{\|\psi_{H}\|_{\hildisnorm} = 1} \| (\interpol \nohom{u}^{\prime \prime}, \psi_H)_{\hilnormdiff} \|_{L^{\infty}[0,T]}\\
            E_{H,3} & = \|(\operatorname{I} - \interpol)\nohom{u} \|_{\Linf{\solSpac}} + \|(\operatorname{I} - \interpol) \nohom{u}^{\prime} \|_{\Linf{\solSpac}} +\|(\operatorname{I} - \interpol) \nohom{u}^{\prime \prime} \|_{\Linf{\hilSpac}}  \\
            E_{H,4} & = \|\ritzH \nonlin(\nohom{u}, \nohom{u}^{\prime}) - \nonlindis( \ritzV \nohom{u}, \interpol \nohom{u}^{\prime}) \|_{\Linf{\hildisSpac}}   + \|\ritzH\assdampOp\nohom{u}^{\prime} - \assdisdampOp\interpol \nohom{u}^{\prime}\|_{\Linf{\hildisSpac}} . \\
        \end{aligned}
    \end{equation}
            \end{subequations}
\end{Satz}
In order to prove the Theorem, we first define the remainder term, which is an essential part of the error.
\begin{Definition}
    \label{errorbounds::remainder}
    In the following let $z=[u,v]^T \in \hilSpacblock$.
    \begin{enumerate}[label=\textnormal{(\arabic*)}]
        \item[(i) ]\label{errorbounds::remainer::1} The remainder of the linear monotone operator is given by
        \begin{align*}
            \remainderdiffblock&:\refSpacblock \cap  D(\diffblockrestrict)\to \soldisSpacblock,                                                                        \\   \remainderdiffblock z  &:= \ritzblock \diffblockrestrict z  - \diffdisblock \interpolblock z 
            = \left[ \begin{array}{rrr}
                \ritzV v- \interpol v               &                             \\
                \ritzH (\assOp u + \assdampOp v ) - (\assdisOp\ritzV u + & \hspace*{-0.3cm}\assdisdampOp \interpol v ) \\
            \end{array}\right]. 
        \end{align*}
        \item[(ii) ]\label{errorbounds::remainder::2}The remainder of the Lipschitz continuous nonlinearity is given by
        \begin{align*}
            \remainderlipschitz (z): \refSpacblock \to \soldisSpacblock:\qquad	\remainderlipschitz (z) :&= \ritzblock \nonlinblock (z) -\nonlinblockdis(\interpolblock (z))
            \\ & = \left[ \begin{array}{rrr}
                0                    &                             \\
                \ritzH \nonlin(u,v)- \hspace*{-0.3cm} & \nonlindis(\ritzV u, \interpol v) + \ritzH f - f_H \\
            \end{array}\right].
        \end{align*}
    \end{enumerate}
\end{Definition}
The next lemma provides an error bound for the remainder term. The proof can be found \cite[Lemma 4.7]{Hipp2019}. Note that the reference uses  additional lift operators $Q_h^V, Q_h^H$, which are here the identity.
\begin{Lemma}
    \label{errorbounds::remainderestimate}
    For  $z = [u,v]^T \in ( \refSpac \times \refSpac) \cap D(\diffblockrestrict)$ we have
    \begin{align*}
        \|\remainderdiffblock z \|_{\soldisSpacblock} &\leq \constremainder \Bigl(
        \max_{\|\phi_{H}\|_{\energytdisnorm} = 1} \| (\interpol \nohom{u}, \phi_H)_{\energynormdiff}\|_{L^{\infty}[0,T]}
        + \max_{\|\psi_{H}\|_{\hildisnorm} = 1} \| (\interpol \nohom{u}, \psi_H)_{\hilnormdiff} \|_{L^{\infty}[0,T]}  \\
        &\qquad \qquad +\max_{\|\phi_{H}\|_{\energytdisnorm} = 1} \| (\interpol \nohom{v}, \phi_H)_{\energynormdiff}\|_{L^{\infty}[0,T]} 
         + \|\ritzH\assdampOp\nohom{v} - \assdisdampOp\interpol \nohom{v}\|_{\Linf{\hildisSpac}}
        \\
        &\qquad \qquad 
        +  \|(\operatorname{I} - \interpol)\nohom{u} \|_{\Linf{\solSpac}}  +  \|(\operatorname{I} - \interpol) \nohom{v}  \|_{\Linf{\solSpac}} \Bigr) ,
    \end{align*}
    where $\constremainder$ only depends  constants, which appear in  Section $\ref{subsec:wellposed}$ and  Assumptions $\ref{Spacedis::assumptions}$, $\ref{Spacedis::Operators}$.
\end{Lemma}
In order to provide error bounds, we start by defining  defects
$\defectstagesIMEXRH{n+ \frac{1}{2}}$ and $\defectIMEXRtotalH{n+1}$ of the IMEX scheme  \eqref{errorbounds::IMEXR} via
    \begin{subequations}
        \label{fulldis::IMEX:exact}
        \begin{alignat}{2}
            \solvechalbntildeH &= \solvecntildeH + \frac{\tau}{2}\diffdisblock \solvechalbntildeH  + \frac{\tau}{2}\nonlinblockntildeH + \frac{\tau}{2} \defectstagesIMEXRH{n+ \frac{1}{2}}, \label{fulldis::IMEX:exact1}  \\
            \solvecnptildeH &= \solvecntildeH +  \tau \diffdisblock \solvechalbntildeH + \tau \nonlinblocknhalbetildeH + \tau \defectIMEXRtotalH{n+1}. \label{fulldis::IMEX:exact2}
        \end{alignat}
    \end{subequations}	
            Here, we used the short notation
    \begin{equation}
        \label{fulldis::notation1}
        \begin{aligned}
            \solvecntildedis = \interpolblock\solvec(t_n), \qquad \nonlinblockntildeH = \nonlinblockdis(\tilde{\solvec}_H^n). 
        \end{aligned}
             \end{equation}
    Bounds on the defects are given in the following Lemma in terms of the constants
\begin{align}
    \constderivative{j} & = \norm{\solvec^{(j)} }_{\L{\infty}{\hilSpacblock}}, &
    \constSderivate & = \norm{\diffblock \solvec^{\prime\prime}}_{\L{\infty}{\hilSpacblock}},& 
\end{align}
and $E_H$ defined in \eqref{eq:EHdef}.
\begin{Lemma} \label{lem:defects}
    Let  Assumptions of Theorem~$\ref{fulldis::TheoremIMEX}$ hold true.
    Then the defects defined in \eqref{fulldis::IMEX:exact} satisfy
    \begin{subequations}
        \begin{align}
            \tauhalb \normXH{\defectstagesIMEXRH{n + \frac{1}{2}}} & \leq  \constdisX  \bigl( \constderivative2
            + \constlipschitzdis  \constderivative{1}\bigr) \frac{\tau^2}{4} + \tauhalb \constdefects  E_H, \\
            \tau \normXH{\defectIMEXRtotalH{n+1}} & \leq \constdisX \constderivative{3}  \,\frac{\tau^3}{24} 
            + \tau \constdefects  E_H, \\ 
            \tauhalb\normXH{\defectstagesIMEXRdiffH{n+\frac{1}{2}} } & \leq \constdisX  ( \constderivative{3} + \constSderivate  ) \frac{\tau^3}{2} + \tau \constdefects  E_H , 
        \end{align}
    \end{subequations}
    where $\defectstagesIMEXRdiffH{n+\frac{1}{2}} = \defectstagesIMEXRH{ n  + \frac{1}{2} } - \defectstagesIMEXRH{ n  - \frac{1}{2}}$ and $\constdefects = \max\{1, \constdisH, \constremainder\}$.
\end{Lemma}
A key idea of the proof is to interpret the IMEX scheme \eqref{errorbounds::IMEXR} as a perturbation of the implicit midpoint rule. 
\begin{proof}
\textit{(a) } 
Set $ \nonlinblockntilde = \nonlinblock\bigl( \solvec(t_n)\bigr)$ and we get 
\begin{align}
        \tauhalb\defectstagesIMEXRH{n + \frac{1}{2}}
        & = \solvechalbntildeH - \solvecntildeH - \frac{\tau}{2} \bigl(\diffdisblock \solvechalbntildeH  + \nonlinblockntildeH \bigr) 
        \nonumber \\
     &= 
     \tauhalb\ritzblock \bigl(\defectstagesCN{n + \frac12}   
         + \nonlinblockhalbetilde{n} - \nonlinblocktilde{n}   \bigr)
        + \tauhalb \defectspaceMPstage{n + \frac12}  ,
        \label{eq:DIMEX}
\end{align}
where 
\begin{align*}
     \tauhalb\defectstagesCN{n + \frac12}  = \solvechalbtilde{n} -\solvecitilde{n} - \frac{\tau}{2}(\diffblock \solvechalbtilde{n}  + \nonlinblocktilde{n}) &= \int_{t_n}^{t_{n + \frac{1}{2}}}  \solvec^{\prime}(s) \d{s} - \tauhalb \solvec^{\prime}(t_{n + \frac12})
\end{align*}
is the projected defect arising in the first stage of the implicit midpoint rule \eqref{Setting::FoF:1} applied to \eqref{Setting::fof} and 
\begin{equation}  \label{eq:defectspaceMP}
    \begin{aligned}
        \tauhalb \defectspaceMPstage{n + \frac12} &=  (\interpolblock - \ritzblock ) (\solvecitilde{n + \frac12} - \solvecitilde{n})   + \tauhalb  \bigl(\remainderdiffblock\solvecitilde{n+\frac{1}{2}}+ \remainderlipschitz(\solvecitilde{n }) \bigr) \\ 
         &=  \tauhalb\Bigl( (\interpolblock - \ritzblock )   \int_{0}^{1}  \solvec^\prime(t_n + \tauhalb s) \d{s}   + \remainderdiffblock\solvecitilde{n+\frac{1}{2}}+ \remainderlipschitz(\solvecitilde{n }) \Bigr) . 
    \end{aligned}
\end{equation}
Since $\tauhalb\defectstagesCN{n + \frac12}$ is the error of the rectangular rule applied to $\solvec^\prime$ on $[t_n,t_{n+\frac12}]$, we can bound it by 
\begin{equation}\label{eq:defectsMP}
    \tauhalb\normX{\defectstagesCN{n + \frac12}}
     \leq \frac{\tau^2}{8} \norm{\solvec^{\prime\prime} }_{\L{\infty}{\hilSpacblock}} .
\end{equation}
Thus, \eqref{eq:DIMEX} and the Lipschitz continuity (Assumption~\ref{Setting::assumption})  yield
\begin{equation*}
        \tauhalb \normXH{\defectstagesIMEXRH{n + \frac{1}{2}}}
        \leq \constdisX  \frac{\tau^2}{4}\bigl( \norm{\solvec^{\prime\prime} }_{\L{\infty}{\hilSpacblock}}
                    + \constlipschitzdis  \norm{\solvec^{\prime} }_{\L{\infty}{\hilSpacblock}}\bigr)
                    + \tauhalb\normXH{\defectspaceMPstage{n + \frac12}},
\end{equation*}
Moreover, by Lemma \ref{errorbounds::remainderestimate} and \cite[Theorem 4.8]{Hipp2019} we get 
\begin{equation}
    \label{eq:spaceerror_estimate}
    \begin{aligned}
            \normXH{\defectspaceMPstage{n + \frac12} } &\leq \tauhalb \normXH{\int_{0}^{1}(\interpolblock - \ritzblock )\solvec^\prime(t_n + \tauhalb s) \d{s} +   \remainderdiffblock(\solvecitilde{n+\frac{1}{2}}) + \remainderlipschitz(\solvecitilde{n})} \\ 
            &\leq \tauhalb \Bigl(  \constdisH\|( \identity - \interpol ) u^{\prime\prime}\|_{\L{\infty}{\hilSpac}}  + \max_{\|\psi_{H}\|_{\hildisnorm} = 1} \| (\interpol \nohom{u}^{\prime \prime}, \psi_H)_{\hilnormdiff} \|_{L^{\infty}[0,T]} \\ & \qquad  \qquad +  \|\remainderdiffblock(\solvecitilde{n+\frac{1}{2}})\|_{\soldisSpacblock} + \|\ritzH \nonlin(\nohom{u}, \nohom{u}^{\prime}) - \nonlindis( \ritzV \nohom{u}, \interpol \nohom{u}^{\prime}) \|_{\Linf{\hildisSpac}} \Bigl) \\ 
            &\leq \tauhalb \max\{1, \constdisH, \constremainder\} E_H.
    \end{aligned}
\end{equation}
\textit{(b) } 
Next we bound the defect $\defectIMEXRtotalH{i}$ defined in \eqref{fulldis::IMEX:exact2}. We obtain
\begin{align}
    \tau\defectIMEXRH{n + 1}
    & = \solvecitildeH{n+1} - \solvecntildeH - \frac{\tau}{2} \bigl(\diffdisblock \solvechalbntildeH  + \nonlinblocknhalbetildeH \bigr) 
    \nonumber \\
    &= 
    \tau \ritzblock \defectCN{n + 1}   
    + \tau \defectspaceMP{n + 1} ,
    \label{eq:dIMEX}
\end{align}
where 
\begin{align*}
     \tau\defectCN{n + 1}  = \solvecitilde{n+1} -\solvecitilde{n} - \tau(\diffblock \solvechalbtilde{n}  + \nonlinblocktilde{n + \frac12}) &= \int_{t_n}^{t_{n + 1}}  \solvec^{\prime}(s) \d{s} - \tau \solvec^{\prime}(t_{n + \frac12})
\end{align*}
is the projected defect arising in the second stage of the implicit midpoint rule \eqref{Setting::FoF:2} applied to \eqref{Setting::fof} and 
\begin{equation*}
    \begin{aligned}
        \tau \defectspaceMP{n + 1} =  (\interpolblock - \ritzblock ) (\solvecitilde{n + 1} - \solvecitilde{n})   + \tau  \bigl(\remainderdiffblock(\solvecitilde{n+\frac{1}{2}}) + \remainderlipschitz(\solvecitilde{n +  \frac12})\bigr).
    \end{aligned}
\end{equation*}
Analogously to \eqref{eq:spaceerror_estimate} we conclude
\begin{equation}\label{eq:spaceerror_estimate2}
    \begin{aligned}
              \tau \normXH{\defectspaceMP{n + 1}}  \leq \tau \max\{1, \constdisH, \constremainder\} E_H.
    \end{aligned}
\end{equation}
The quadrature error of the midpoint rule satisfies
\begin{equation*}
    \tau\normX{\defectCN{n + 1}}
    \leq \frac{\tau^3}{24} \norm{\solvec^{\prime\prime\prime} }_{\L{\infty}{\hilSpacblock}} .
\end{equation*}
This, \eqref{eq:dIMEX}, and \eqref{eq:spaceerror_estimate2} yield
\begin{equation}
    \begin{aligned}
            \tau \normXH{\defectIMEXRH{n + \frac{1}{2}}}
            \leq \constdisX  \,\frac{\tau^3}{24} \norm{\solvec^{\prime\prime\prime} }_{\L{\infty}{\hilSpacblock}}
                        + \tau \max\{1, \constdisH, \constremainder\} E_H, 
    \end{aligned}
\end{equation}
(c)  Finally, we consider $\defectstagesIMEXRdiffH{n+\frac{1}{2}}$. 
From \eqref{eq:DIMEX} we get 
\begin{equation}
\begin{aligned}
    \label{eq:defectdiffMPH}
    \tauhalb\defectstagesIMEXRdiffH{n + \frac12} =  \tauhalb \ritzblock (\defectstagesCNdiff{n+\frac{1}{2}} &+ (\nonlinblockhalbetilde{n} - \nonlinblocktilde{n} ) -(\nonlinblockhalbeminustilde{n} - \nonlinblockminustilde{n}) ) \\ &+ (\defectspaceMPstage{n+ \frac12} - \defectspaceMPstage{n - \frac12})  . \\ 
    %
\end{aligned}
\end{equation}
Analogously to \eqref{eq:defectspaceMP}, the defect $ \defectstagesCNdiff{n + \frac12}$ is bounded by
\begin{equation}
    \label{eq:defectdiffMP}
    \begin{aligned}
        \tauhalb\normX{\defectstagesCNdiff{n + \frac12}}
        \leq \frac{\tau^3}{8} \norm{\solvec^{\prime\prime\prime} }_{\L{\infty}{\hilSpacblock}} .
    \end{aligned}
\end{equation}
Using \eqref{Setting::fof} and Taylor expansion, one can easily verify that
\begin{equation*}
    \begin{aligned}
        \tauhalb\normX{ (\nonlinblockhalbetilde{n} - \nonlinblocktilde{n} ) -(\nonlinblockhalbeminustilde{n} -  \nonlinblockminustilde{n}) }
        &\leq  \Bigl(\tauhalb\Bigr)^2 \tau  \Bigl(\norm{\solvec^{\prime \prime \prime} }_{\L{\infty}{\hilSpacblock}} + \norm{  \diffblock \solvec^{\prime \prime} }_{\L{\infty}{\hilSpacblock}} \Bigr).
    \end{aligned}
    \end{equation*}
From this, the space error estimate  \eqref{eq:spaceerror_estimate} and the relation \eqref{eq:defectdiffMPH}, it follows that:
\begin{equation*}
    \begin{aligned}
        \tauhalb\normXH{\defectstagesIMEXRdiffH{i+\frac{1}{2}} }
        &\leq  \constdisX \frac{\tau^3}{8} ( 3\norm{\solvec^{\prime \prime \prime} }_{\L{\infty}{\hilSpacblock}} + 2\norm{  \diffblock \solvec^{\prime \prime} }_{\L{\infty}{\hilSpacblock}}  )  + \tauhalb (\|\defectspaceMPstage{n +\frac12 }  \|_{\soldisSpacblock} +  \|\defectspaceMPstage{n - \frac12}  \|_{\soldisSpacblock }) \\
        &\leq  \constdisX \frac{\tau^3}{2} ( \norm{\solvec^{\prime \prime \prime} }_{\L{\infty}{\hilSpacblock}} + \norm{  \diffblock \solvec^{\prime \prime} }_{\L{\infty}{\hilSpacblock}}  )+ \tau \max\{1, \constdisH, \constremainder\}  E_H 
    \end{aligned}
   \end{equation*}
    This completes the proof.
\end{proof}
Next we introduce the operators
\begin{equation}
    \label{errorbounds::ROperators}
    \begin{aligned}
        \resplusminusOpdis \coloneqq \identity \pm \tauhalb \diffdisblock: \soldisSpacblock \to \soldisSpacblock.
    \end{aligned}
\end{equation}
The following properties of $\resplusminusOpdis$ can be shown similarly to the continuous case \cite[Lemma 2.4]{Leibold2021Hochbruck}:
\begin{Lemma} \label{errorbounds::propertiesR} 
Let $\constqmdisS$ be as defined in Lemma~$\ref{Settig::lemma:fo_qmono}$.
Then, if $\tau \constqmdisS < 2$, the following assertions hold true:
\begin{enumerate}
    \item[ a) ] $\resminusOpdis$ is invertible with $\normOperatorXH{\resminusOpinvdis} \leq 1$ 
    \item[ b) ] $\resOpdis\coloneqq \resplusOpdis \resminusOpinvdis$ satisfying $\normOperatorXH{\resOpdis} \leq \e^{\tau \constqmdisS}$.
\end{enumerate}
\end{Lemma}

The main theorem can now be proved as follows.
\begin{proof}[Proof of Theorem $\ref{fulldis::TheoremIMEX}$]
    $\,$ We start by splitting the error into two parts
    \begin{equation}
        \label{fulldis::spliterror}
        \begin{aligned}
            \solvecndis - \solvec(t_n) = \errornH + \errorH,
        \end{aligned}
    \end{equation}
    where 
    \begin{equation*}
        \begin{aligned}
            \errornH &= \solvecndis - \interpolblock \solvec(t_n), \\
            \errorH      &= \interpolblock \solvec(t_n) - \solvec(t_n).
        \end{aligned}
    \end{equation*}
    As in \cite[equation (27)]{Leibold2021Hochbruck} we have
    \begin{equation}
        \label{fulldis::errorboundprojection}
        \begin{aligned}
            \normX{\errorH } \leq C(\abstractspaceError{2} + \abstractspaceError{3}).
        \end{aligned}
    \end{equation}
    Hence, it remains to bound $\errornH$.

            (a) \textit{Error recursion.} By subtracting \eqref{fulldis::IMEX:exact} from \eqref{errorbounds::IMEXR} we obtain the error recursion
    \begin{subequations}
        \label{fulldis::IMEX:error}
        \begin{align}
            \errornhalbH &= \errornH  + \frac{\tau}{2}\diffdisblock \errornhalbH  + \frac{\tau}{2}\errornonlinearnH -  \frac{\tau}{2} \defectstagesIMEXRH{n+ \frac{1}{2}} \label{fulldis::IMEX:error1}, \\
            \errornpH &= \errornH  + \tau \diffdisblock \errornhalbH + \tau \errornonlinearnhalbeH - \tau \defectIMEXRH{n+1}, \label{fulldis::IMEX:error2}
        \end{align}
    \end{subequations}
    where $\errornonlinearnH = \nonlinblockntildeH - \widetilde{F}_H^{n-1}$.
    We rewrite \eqref{fulldis::IMEX:error1} using the operators defined in \eqref{errorbounds::ROperators} as
    \begin{equation}
        \label{fulldis::IMEX:errorhalbwithR}
        \begin{aligned}
            \errornhalbH = \resminusOpinvdis(\errornH + \tauhalb (\errornonlinearnH -   \defectstagesIMEXRH{n+ \frac{1}{2}}) ).
        \end{aligned}
    \end{equation}
    Plugging \eqref{fulldis::IMEX:errorhalbwithR} into \eqref{fulldis::IMEX:error2} and using $\tau \diffdisblock = \resplusOpdis - \resminusOpdis$ yields
    \begin{align}
        \errornpH &=  \resOpdis \errornH  + \tauhalb ( \resOpdis - \operatorname*{I}) (\errornonlinearnH -   \defectstagesIMEXRH{n+ \frac{1}{2}}) + \tau \errornonlinearnhalbeH - \tau \defectIMEXRH{n+1}. \label{fulldis::IMEX:recursion}
    \end{align}

    (b) \textit{Stability.} 
    Solving the error recursion \eqref{fulldis::IMEX:recursion} and using Lemma \ref{errorbounds::propertiesR} yields
    \begin{equation*}
        \begin{aligned}
            \normXH{\errornH} 
            &=  && \normXH{\resOpdis^n e_H^0}  + \tauhalb \suminullnminuseins \e^{\tau \constqmdisS (n - 1 -i)} \normOperatorXH{ \operatorname*{I}  - \resOpdis  } \normXH{\errornonlineariH}    \\
            & &&+\tauhalb   \normXH{  \suminullnminuseins (\resOpdis^{n-i} - \resOpdis^{n-i-1} )  \defectstagesIMEXRH{i+ \frac{1}{2}}} \\
            & &&+ \tau \suminullnminuseins \e^{\tau \constqmdisS (n - 1 -i)}  \normXH{\errornonlinearihalbeH}
             + \tau \suminullnminuseins  \e^{\tau \constqmdisS (n - 1 -i)}\normXH{\defectIMEXRH{i+1}}
             .\label{fulldis::nachVCF} 
        \end{aligned}
    \end{equation*}
    Multiplying by $\e^{-\tau \constqmdisS n} $ on both sides, using the Lipschitz continuity and summation by parts for $\defectstagesIMEXRH{n}$ terms results in
    \begin{equation}
        \label{fulldis::bringtootherside}
        \begin{aligned}
            \e^{-\tau \constqmdisS n}\normXH{\errornH} \leq \quad & \normXH{e_H^0 } +\constlipschitzdis \frac{\tau}{2} (  \e^{-\tau \constqmdisS } + 1)\suminullnminuseins \e^{- \tau \constqmdisS i}  \normXH{\erroriH}    \\
            + & \, \tauhalb  \bigl( \normXH{   \defectstagesIMEXRH{1/2} } +   \normXH{\defectstagesIMEXRH{n - 1/2} } 
            +    \sumi{1}{n-1}\normXH{  \defectstagesIMEXRdiffH{i+\frac{1}{2}}} \bigr) \\
            + & \,\tau 
            \suminullnminuseins  \normXH{ 
            \defectIMEXRH{i+1}} 
            + \constlipschitzdis \tau\suminullnminuseins \e^{- \tau \constqmdisS i}  \normXH{\errorihalbH} .
        \end{aligned}
    \end{equation}
    Next, we need to estimate the error of the half step in terms of the full step. By \eqref{fulldis::IMEX:errorhalbwithR}, the Lipschitz continuity and Lemma \ref{errorbounds::propertiesR}, we derive
    \begin{equation*}
        \begin{aligned}
            \normXH{\errornhalbH}  &\leq (1+ \tauhalb \constlipschitzdis )\normXH{\errornH}  +    \tauhalb \normXH{\defectstagesIMEXRH{n+ \frac{1}{2}}}. 
        \end{aligned}
    \end{equation*}
    Using this estimate in \eqref{fulldis::bringtootherside} and $ \e^{- \tau \constqmdisS i}  \leq 1$ yields
    \begin{equation*}
        \begin{aligned}
            \e^{-\tau \constqmdisS n}\normXH{\errornH} \leq   \constlipschitzdis & \tau (2+ \tauhalb \constlipschitzdis ) \suminullnminuseins \e^{- \tau \constqmdisS i}  \normXH{\erroriH} \\ & + \normXH{e_H^0 } + \abstractfullError{1} + \abstractfullError{2} + \abstractfullError{3},
        \end{aligned}
    \end{equation*}
    where 
    \begin{equation}
        \label{fulldis::abstracttime}
        \begin{aligned}
        \abstractfullError{1} &= \constlipschitzdis \tau \suminullnminuseins \tauhalb \normXH{\defectstagesIMEXRH{i + \frac{1}{2}}}  \\
        \abstractfullError{2} &= \tauhalb \Bigl(\normXH{  \defectstagesIMEXRH{1/2}} +  \normXH{ \defectstagesIMEXRH{n - 1/2}} + \sumi{1}{n-1}  \normXH{ \defectstagesIMEXRdiffH{i + \frac{1}{2}}} \Bigr)\\
        \abstractfullError{3} &=  \tau \suminullnminuseins  \normXH{\defectIMEXRH{i+1}}.
        \end{aligned}
    \end{equation}
    By the Grönwall's inequality and $\normXH{e_H^0} \leq \abstractspaceError{1}$, we get
    \begin{equation}
        \label{fulldis::IMEX:afterGronwall}
        \begin{aligned}
            \normXH{\errornH} \leq C(  \abstractspaceError{1} + \abstractfullError{1} + \abstractfullError{2} + \abstractfullError{3}) \e^{\bigl(  \constlipschitzdis (2+ \tauhalb \constlipschitzdis )t_n + \constqmdisS t_n \bigr)}.
        \end{aligned}
    \end{equation}
    Overall, we get by Lemma \ref{lem:defects}, \eqref{fulldis::errorboundprojection}, \eqref{fulldis::IMEX:afterGronwall}, and 
    \begin{equation*}
        \begin{aligned}
            \| \nohomdisn{u} - \nohom{u}(t_n)\|_{\energytnorm} + \| \nohomdisn{v} - \nohom{u}^{\prime}(t_n) \|_{\hilnorm} \leq \normX{\solvecnH - \solvec(t_n)}, 
        \end{aligned}
    \end{equation*}
    the error bound \eqref{fulldis::IMEX:errorboundtotal}.
    \end{proof}
    Under additional assumptions, the existence of $\nohomdisn{u}$ and $\nohomdisn{v}$ on $[0,T]$ can be proved analogously to \cite[Corollary 3.5]{Leibold2021Hochbruck} by the following Theorem. 
    \begin{Satz}
        \label{fulldis::IMEX:Theorem:bound}
        Let Assumption of Theorem $\ref{fulldis::TheoremIMEX}$ be fulfilled. 
        Further we assume that the spatial error $E_H$
                    defined in \eqref{errorbounds::abstractspaceError} satisfies $E_H \to 0$ for $H \to 0$ and set\[\rho \coloneqq 3  ( \constdisV \|\nohom{u}\|_{\L{\infty}{\solSpac}} + \constdisH\|\nohom{u}^\prime \|_{\L{\infty}{\hilSpac}} ) + \constinterpol\|\nohom{u}^{\prime}\|_{\L{\infty}{\refSpac}}.\] Then, there exist $\tau^\star,H^\star>0$ such that for all $h < H^\star$, $\tau < \tau^\star$ the approximation $ [\nohomdisn{u} \, \nohomdisn{v}]$ obtained by the fully discrete IMEX scheme satisfies
        \[ \max_{t_n \leq T} \| \nohomdisn{u}\|_{\energytdisnorm} \leq \rho \quad \text{  and  } \quad \max_{t_n \leq T} \| \nohomdisn{v}\|_{\hildisnorm} \leq \rho  \]
        and therefore the error bound 	\eqref{fulldis::IMEX:errorboundtotal}
        is valid with $\rho_H = \rho$.
    \end{Satz}
    \begin{Bemerkung} \phantom{2} 
        \label{fulldis::remark}
        \begin{enumerate}
            \item[(i) ]  Theorem \ref{IMEX::Theorem} can be proved in a similar way as above by replacing the discrete expressions with the continuous ones, since the assumptions are almost identical in both cases. The error terms in space vanish in the semidiscrete case.  
            \item[(ii) ] Analogously, one can obtain that  Theorems \ref{fulldis::TheoremIMEX} and \ref{fulldis::IMEX:Theorem:bound} are also valid for the implicit midpoint rule if $\tau < \min\{\frac{1}{\constqmdisS}, \frac{2}{\constlipschitzdis}\}$ with $M = \constqmdisS + \frac{\constlipschitzdis}{1 - \tauhalb\constlipschitzdis}$. 
        \end{enumerate}
    \end{Bemerkung}
\section{Application: Semilinear damped wave equations with highly oscillatory coefficients \label{Application}}
In this section, we illustrate the theory of the previous sections using a specific nonlinear model problem.
In particular, we are interested in wave problems with nonlinear damping terms and a coefficient which is oscillating fast in space. Using G-convergence results we derive a homogenized version, for which we apply the heterogeneous multiscale method.
\subsection{Notations}
Let $\Omega \subset \R^d, d \leq 3$, be a bounded domain with a Lipschitz boundary $\partial \Omega$. With $W^{s,p}(\Omega)$ we denote the standard Sobolev space and with $W_0^{s,p}(\Omega)$ the space consisting of functions of $W^{s,p}(\Omega)$ which are zero on $\partial \Omega$ in the trace sense. If $p=2$ we use the notation $\Hspace{s}$ and $\Hzerospace{s}$. 
In addition, let  \[ H^{s}_{\#}(Y) \coloneqq \Bigl\{  w \in H^s_{per} \mid \int_{Y} w = 0 \Bigr\}, \]
 where  $ H^s_{per}(Y)$ is the closure of $ C^{\infty}_{per}(Y)$ 
 in the $H^s$-norm on a hypercube $Y \subset \Omega$. In the following, we denote $Y_{\sigma}$ as the hypercube with side length $\sigma$ and $Y$ as the unit hypercube.
 Further, there exists a constant $\constpoincare>0$, such that the Poincare inequality is fulfilled  $w\in H^1_0(\Omega)$.
\[ \| w \|_{\Lspace{2}} \leq \constpoincare \|\nabla w\|_{\Lspace{2}}. \] For a Banach space $B$ we denote by $\L{p}{B}$ the Bochner space of functions $u: [0,T] \to B$ equipped with the norm \[\|u\|_{\L{p}{B}} \coloneqq \Bigl(\int_0^T \|u(t)\|^p_B \d{t}\Bigr)^{\frac{1}{p}}. \]
In the following, we set \[ \solSpac = \Hzerospace{1} \quad \text{ and } \quad \hilSpac = \Lspace{2} .\] 
\subsection{Model problem and G-convergence}\label{sec:model}
We consider the semilinear wave equation with damping:
\begin{equation}
    \label{model::modelproblem}
    \left\{\begin{aligned}
        \partial_{tt} u^{\veps} - \nabla\cdot (\Coeffeps(x)\nabla u^{\veps})  - \nabla\cdot ( \beta(x) \nabla \partial_tu^{\veps})   & =   \nonlin(x,\partial_{t}u^{\veps})+ f             &                    & \text{in } (0,T) \times\Omega,           \\  
        u^{\epsilon}     & = 0             &                    & \text{on } (0,T) \times \partial \Omega, \\
        \hspace*{1.45cm} u^{\veps}(0) = \nohomnull{u} \hspace*{0.56cm}  ,\hspace*{0.56cm} 	\partial_{t}u^{\veps}(0)                                          & = \nohomnull{v}
        &                 & \text{in } \Omega.
    \end{aligned}\right.
\end{equation}
 $\partial_{t} u^{\veps}$ and  $\partial_{tt} u^{\veps}$ are time derivatives of $u$, which we denoted in Section \ref{generalSet} and \ref{fulldis} as $u^\prime$ and $u^{\prime\prime}$ , respectively. However, since the index $\veps$ occurs here, we have changed the notation for better readability.
For the coefficient $\Coeffeps$, $\beta$ and the nonlinear operator  $\nonlineps$ we assume the following
\begin{Annahme}
    Let $\constcoercive,\constbounded, \constBeta, \constbeta$ be positive constants.
    \label{model::assumption}
    \begin{enumerate}[label=\textnormal{(\arabic*)}]
        \item[(i) ] $ \Coeffeps \in (W^{1,\infty}(\Omega))_{sym}^{d\times d}$ and 
        $ \constbounded |\xi|^2 \geq  \xi \cdot \Coeffeps(x) \xi > \constcoercive \,|\xi|^2$, a.e. in $\Omega$ for all $\xi \in \R^n$, $\xi \neq 0$. 
        \item[(ii) ]  $ \beta \in C(\overline{\Omega}) $ and $\constbeta< \beta(x) \leq \constBeta$.
        \item[(iii) ] $ \nonlin(\cdot, \eta)$ is Lebesgue measurable for all $\eta  \in \R$ and 
        \[\nonlin(x,\cdot): \L{2}{\hilSpac} \to \L{2}{\hilSpac}  \] maps bounded sets to bounded sets for almost all $x \in \Omega$. 
        \item[(iv) ] $\| \nonlin(\cdot,v)-\nonlin(\cdot,w)\|_{\hilSpac} \leq \constlipschitz \|v-w\|_{\hilSpac}$ for all $v,w \in \hilSpac$ with $\|v\|_{\hilSpac},\|w\|_{\hilSpac} \leq \rho$. 
        \item[(v) ] For all sequences $\{w^{\veps}\}$ with $w^{\veps} \rightarrow w$ in $\L{2}{\hilSpac}$ for $\veps \rightarrow 0$ there exists a subsequence which we still denote by $\{w^{\veps}\}$ with \[\nonlin( \cdot,w^{\veps}) \rightharpoonup \nonlin(\cdot, w)\] in  $\L{2}{\hilSpac}$.
        \item[(vi) ] 	
        \begin{align*}
            f \in W_{loc}^{1,1}([0,\infty);H^2(\Omega) ) \quad \text{ and } \quad  \nohomnull{u},\nohomnull{v} \in H^2(\Omega) \cap \solSpac    \asseqqed
        \end{align*}
    \end{enumerate}
\end{Annahme}
The weak formulation of \eqref{model::modelproblem} reads as follows 
\begin{equation}
    \label{model::weakformulation}
    \begin{aligned}
        \bigl(\partial_{tt} {u}^{\veps}, \Phi \bigr)_\hilnorm + \bilinear{}^{\veps}({u}^{\veps}, \Phi) + \bilinearbeta{}( \nabla \partial_t {u}^{\veps}, \nabla \Phi )  = ( \nonlin(\cdot,\partial_t {u}^{\veps}), \Phi )_\hilnorm  +  \bigl( f,\Phi\bigr)_\hilnorm,
    \end{aligned}
\end{equation}
where 
\begin{equation}
    \label{modelproblem:bilinearforms}
    \begin{aligned}
         \bilinear{}^{\veps}({u}^{\veps}, \Phi) & \coloneqq \bigl(\Coeffeps\nabla{u}^{\veps}, \nabla\Phi\bigr)_{\hilSpac}, \\
         \bilinearbeta{}(u^{\veps},\Phi) &\coloneqq \bigl(\beta \nabla{u}^{\veps}, \nabla\Phi\bigr)_{\hilSpac}.
    \end{aligned}
\end{equation}
for all $\Phi \in \solSpac$.
Next we check whether \eqref{model::weakformulation} fits into the abstract setting.  
\begin{Lemma}
    \label{model::Lemma:assumption}
    Let Assumption $\ref{model::assumption}$ be satisfied. Then the weak formulation $\eqref{model::weakformulation}$ is  of the form $\eqref{Setting::weak}$. Further Assumption $\ref{Setting::assumption}$ is satisfied with $\constsp = 1$ and $\constqm = 0$. 
\end{Lemma}
\begin{proof}
     $\,$ This is a direct consequence of  Assumption \ref{model::assumption} and \eqref{modelproblem:bilinearforms}. 
\end{proof}
By Lemma \ref{model::Lemma:assumption} and Corollary \ref{Setting::Corollary:wellposedness} \eqref{model::modelproblem} is locally wellposed.\\
Next we give an example for which the Assumptions \ref{model::assumption}(iii),(iv) are satisfied.
\begin{Bsp}
    \label{model::Gexample}
    We set
    \begin{equation}
        \begin{aligned}
            \nonlin(x, \eta) = \theta(x)\bigl((|\eta| + \sigma)^{\gamma} - \sigma^{\gamma}\bigr) \sgn(\eta),
        \end{aligned}
    \end{equation}
     $0< \gamma \leq 1$ and $\sigma \geq 0$ where $\theta \in L^{\infty}(\Omega)$. Let $w^{\veps} $ with $ w^{\veps} \rightarrow \nohom{w}$ in $\L{2}{\hilSpac}$. We have to show that, if $\veps \to 0$, 
    \begin{align}
        \label{model::example_convergence}
        \bigl|\int_{(0,T)\times \Omega} (\nonlin(x,w^{\veps})- 		\nonlin(x,\nohom{w}))\varphi\bigr|  \to 0
    \end{align}
    for all $ \varphi \in \L{2}{\hilSpac}$.
    We know that there is a subsequence, which we also denote by $\{w ^{\veps}\}$, such that for almost all $t$
    \begin{equation*}
        \begin{aligned}
            w^{\veps}(t) & \rightarrow \nohom{w}(t) &  & \text{ in } \hilSpac. \\
        \end{aligned}
    \end{equation*}
    According to the Theorem of Yegorov there exist a sequence $\{\mathcal{Q}_{n}\} \subset (0,T) \times \Omega$ with 
    \begin{equation*}
        \begin{aligned}
            \bigl|\bigl((0,T) \times \Omega \bigr) \setminus \mathcal{Q}_{n}\bigr| < \frac{1}{n} \quad \text{ and } \quad w^{\veps} & \rightarrow \nohom{w} &  & \text{ uniform  in  } \mathcal{Q}_{n}. \\
        \end{aligned}
    \end{equation*}
    Further if $\varepsilon$ is small enough
    \begin{equation*}
        \begin{aligned}
            |\nonlin(x,w^{\veps})\varphi|&= |\nohom{\theta}| (( |w^{\veps}|+ \sigma )^\gamma - \sigma^\gamma)|\varphi|  \\ &\leq |\nohom{\theta}| (( |w| +1 + \sigma )^\gamma - \sigma^\gamma)|\varphi| \in L^1(\mathcal{Q}_n).
        \end{aligned}
    \end{equation*}
    By dominated convergence we obtain 
    \begin{align*}
        \int_{\mathcal{Q}_{n}}   \nonlin(x,w^{\veps})\varphi \to \int_{\mathcal{Q}_{n}}  \nonlin(x,\nohom{w})\varphi.
    \end{align*}
    If we now let $n \to \infty$ it follows that  \eqref{model::example_convergence} converges to zero. 
    If $\sigma > 0$, then $\nonlin(x, \cdot) \in C^1(\R)$ for almost all $x \in \Omega$. 
    Then $\nonlin$ is Lipschitz continuous with respect to the $\hilSpac$-norm, since
    \begin{align*}
        \|\nonlin(v) - \nonlin(w)\|_{\hilnorm} & \leq \| \theta\int_0^1 \partial_{\xi}\nonlin(x,v + \xi(v-w)) (v -w)\|_{\hilnorm}                    \\
        & \leq \|\theta\|_{\Lohnespace{\infty}} \| \frac{\gamma}{( |\eta| + \sigma)^{1-\gamma}}\|_{\Lohnespace{\infty}(\R)}\|(v -w)\|_{\hilnorm} \leq  \frac{\gamma \|\theta\|_{\Lohnespace{\infty}}}{ \sigma^{1-\gamma}} \;\|v -w\|_{\hilnorm}.      \bspeqqed
     \end{align*}
\end{Bsp}
By the assumptions we obtain the following convergence result. 
\begin{Satz} \hspace*{-0.2cm} {\normalfont \cite[Theorem 8.3]{Svanstedt2007}}
    \label{model::Theorem:convmodelproblem}
    Consider the sequence of damped hyperbolic problems  \eqref{model::modelproblem} under the  Assumptions $\ref{model::assumption}$.
    Then,
    \begin{equation*}
        \begin{aligned}
            u^{\veps}                               & \rightharpoonup \hom{u}                               &  & \text{ in }  \L{2}{\solSpac}      \\
            \Coeffeps\nabla u^{\veps}      & \rightharpoonup \hom{a}\nabla \hom{u}      &  & \text{ in } \L{2}{\hilSpac} \\
            \nonlin(\partial_{t}u^{\veps}) & \rightharpoonup \nonlin(\partial_{t}\hom{u}) &  & \text{ in } \L{2}{\hilSpac},
        \end{aligned}
    \end{equation*}
    where $\hom{u}$ is the solution of \eqref{model::HomLösunglinear1}.
    \begin{equation}
        \label{model::HomLösunglinear1}
        \left\{\begin{aligned}
            \partial_{tt} \hom{u} - \nabla\cdot ( \hom{a}\nabla \hom{u}) - \nabla\cdot ( \beta \nabla \partial_t\hom{u}) =  \nonlin(\partial_{t}\hom{u}) &+ f  &&\text{in } (0,T) \times\Omega \\ \hspace*{7.28cm}
            \hom{u}                                                                                        & = 0        &                     & \text{on } (0,T) \times \partial \Omega \\
            \hspace*{2.45cm} \hom{u}(0) = \homnull{u} \hspace*{.46cm}  ,\hspace*{0.46cm}  \partial_t\hom{u}(0) & = \homnull{v}
            &            & \text{in }  \Omega.
        \end{aligned}\right.
    \end{equation}
\end{Satz}
\begin{Bemerkung}
 In \cite{Svanstedt2007}, the authors assumed $\beta=1$. However, it is straightforward to extend the proof of the theorem to the more general assumptions, where essentially only constants change. 
\end{Bemerkung}
\begin{Bemerkung}
    \label{model::remark_propertiescoeff}
    It can be seen in the proof that the homogenized coefficient $\hom{a}$ coincides with the homogenized coefficient from an auxiliary parabolic problem. Therefore, we know that the coefficient satisfies Assumption \ref{model::assumption}(i) by parabolic compactness \cite[Theorem 3.1]{Svanstedt1999}. Lemma \ref{model::Lemma:assumption} shows that also the weak formation of the  homogenized equation \eqref{model::HomLösunglinear1} fulfills Assumption \ref{Setting::assumption}. 
\end{Bemerkung}
\subsection{Heterogeneous Multiscale Method}
We formulate the Finite Element Heterogeneous Multiscale Method (FE-HMM). For the macroscopic model we take the homogeneous equation \eqref{model::HomLösunglinear1} which we derived in Section \ref{sec:model}. Based on this model, we want to compute an approximation of $\hom{u}$ by discretization in space. 
We start with the macroscopic problem. We use a partition $\mathcal{T}_H$ into simplicial or quadrilateral elements, where we denote the elements of $\mathcal T_H$ by $K$. As the finite element space we choose \[\disSpace := \{v_{H} \in \solSpac \mid v_{H \mid_{K}} \in \mathbb{P}_p, \forall K \in \mathcal{T}_H\},\]  where $\mathbb{P}_p$ is the space of polynomials of maximal degree $p$.  Further we set $\homdisnull{u} :=I_H\homnull{u}$ and  $\homdisnull{v} := I_H\homnull{v}$, where  $I_H$ is the  interpolation operator to $\disSpace$.
Let $\homdis{u}$ be the solution of the semi-discrete equation
\begin{equation}
    \begin{aligned}
        \label{hmm::diskret}
        \bigl(\partial_{tt} \homdis{u}, \Phi_H \bigr)_{\hilnorm} \! + \! \bilinearhom(\homdis{u}, \Phi_H) \!+ \! \bilinearbeta{}( \partial_t \homdis{u}, \Phi_H ) \!=\! ( \nonlin(\cdot,\partial_t \homdis{u}), \Phi_H )_{\hilnorm}   \!+ \! \bigl( f(t),\Phi_H\bigr)_{\hilnorm}
    \end{aligned}
\end{equation}
for all $\Phi_H \in \disSpace$, with $\homdis{u}(0) = \homdisnull{u}$ and  $\partial_t\homdis{u}(0) = \homdisnull{v}$ with  $\homdisnull{u}, \homdisnull{v} \in \disSpace$.
Theorem \ref{model::Theorem:convmodelproblem} does not provide an explicit representation for the homogenized coefficient $\hom{a}$. However, if we assume more structure on the microscopic scale, i.e., 
\begin{equation}
    \label{model::more_structure}
    \begin{aligned}
        \Coeffeps(x) = a\Bigl(x, \frac{x}{\epsilon} \Bigr) ,
        \end{aligned}
\end{equation}
and  $a(x,y)$ is $Y$-periodic w.r.t $y$, we can derive the same well-known formula for the homogenized coefficient as in the elliptic case \cite[Chapter 1, Eq. 2.20]{Bensoussan1978Papanicolaou} with the help of Remark \ref{model::remark_propertiescoeff} and \cite[Theorem 8.1]{Svanstedt1999}
\[	
    \hom{a}_{ij}(x) = \int_{Y} \sum_{k = 1}^d a_{ik}(x,y)\bigl(\delta_{jk} +				\dfrac{\partial \chi^{j}}{\partial y_k}(x,y)\bigr) \d{y},
    \]
    where $\delta_{jk}$ are the Kronecker delta and $\chi^{j}$ are the solutions of the following cell problems on $Y$  : \\
    Find $\chi^j  \in  H^1_{\#}(Y) $:
    \begin{equation}
        \begin{aligned}
            \int_{Y}a(x, y )   ( e^j + \nabla \chi^j(x,y)) \cdot \nabla z(y) \d{y} = 0
        \end{aligned}
    \end{equation}
    for all $z \in H^{1}_{\#}(Y) $.
 To approximate the resulting bilinear form we use a suitable quadrature formula in the following:
\begin{equation}
    \begin{aligned}
        \bilinearhom(\Phi_H,\Psi_H) = \sum_{K \in \mathcal{T}_H} \int_{K} \hom{a}(x)\nabla \Phi_H \cdot  \nabla \Psi_H 
         \approx \sum_{K \in \mathcal{T}_H} \sum_{j= 1}^{J} \omega^j_{K}\homCoeKj\nabla \Phi_H \cdot  \nabla \Psi_H
        \label{hmm::bilinearform} 
    \end{aligned}
\end{equation}
with the  coefficient
\begin{equation*}
    \begin{aligned}
        \label{hmm::operator}
        \homCoeKj\nabla \Phi_H \cdot  \nabla \Psi_H \coloneqq \hom{a}(x^j_K)\nabla \Phi_H \cdot  \nabla \Psi_H   =   \dfrac{1}{|	\cell{\veps}|}\int_{ \cellKj{\veps}} a\Bigl(x^j_K, \frac{x}{\epsilon} \Bigr) \nabla \phi_j(x) \cdot  	\nabla 	\psi_j(x) \d{x},
                \end{aligned}
\end{equation*}
where $\phi_j,\psi_j$ are the solutions of the shifted cell problems on $\cellKj{\veps}= x^j_K +  Y_\veps$:

Find $\phi_j  - \Phi^{lin}_{H,K,j} \in  H^1_{\#}(\cellKj{\veps}) $:
\begin{equation}
    \label{hmm::cellproblem}
    \begin{aligned}
        \int_{\cellKj{\veps}}a\Bigl(x_K^j, \frac{x}{\epsilon} \Bigr)   \nabla \phi_j(x) \cdot \nabla z(x) \d{x} = 0
    \end{aligned}
\end{equation}
for all $z \in H^{1}_{\#}(\cellKj{\veps}) $, with \[ \Phi^{lin}_{H,K,j}(x) = \Phi_H(x^j_{K}) + (x - x_K^j)\cdot \nabla \Phi_H(x_K^j) .\] $\psi_j$ solves a similar cell problem, using $\Psi_H$ instead of $\Phi_H$. 
\begin{Bemerkung} \label{remark::qf}
    To ensure the optimal convergence rates of the FE-HMM we have to use a suitable quadrature formula $ \{\omega_K^j, x_K^j\}$ of sufficient order.  A detailed discussion including estimates of the quadrature error can be founded \cite[Chp. 4.1]{Abdulle2012Engquist}.
\end{Bemerkung}
    In a next step we discretize the cell problems. Since in applications the finescale parameter $\veps$ may be unknown, we consider as cell domain $\cellKj{\delta}$ with $\delta \geq \veps$.     Further, let $\mathcal{T}_h$ be a partition for each cell domain $\cellKj{\delta} \supseteq\cellKj{\veps}$ consisting of simplicial or quadrilateral elements and a microscopic finite element space consisting of piecewise polynomial functions
    \[\Vmicro(\cellKj{\delta}) := \{v_h \in W(\cellKj{\delta}) \mid v_{h \mid_{\kappa}} \in \mathbb{P}_q, \forall \kappa \in \mathcal{T}_h\},\] where $\mathbb{P}_q$ consist of all polynomials of order $q$ and $W(\cellKj{\delta})$ is a subspace of the Sobolev space $H^{1}(\cellKj{\delta})$. This space describes the coupling between the macroscopical and microscopical problem. Common choices are
    \begin{enumerate}[label=\textnormal{(\arabic*)}]
        \item $\,$ periodic boundary conditions: $W(\cellKj{\delta}):= H^{1}_{\#}(\cellKj{\delta})$
        \item $\,$ Dirichlet boundary conditions: $W(\cellKj{\delta}):= H^{1}_{0}(\cellKj{\delta})$
        \item $\,$ Neumann boundary conditions: $W(\cellKj{\delta}):= H^{1}_{N}(\cellKj{\delta}) = \{u \in \Hspace{1} \mid \int_{\cellKj{\delta}} \nabla u  = \mathbf{0} \} $.
    \end{enumerate}
    We consider the discrete micro problems:
    find $\phihj - \Phi^{lin}_{H,K,j} \in \Vmicro(\cellKj{\delta}) $ :
    \begin{align}
        \label{hmm::cellproblem_dis}
        \int_{\cellKj{\delta}}\Coeffeps(x)  \nabla \phihj(x) \cdot  \nabla z_h(x) \d{x} = 0
    \end{align}
    for all $z_h \in \Vmicro(\cellKj{\delta})$. \\
    Using this, we formulate the discrete form of \eqref{hmm::bilinearform} as
    \begin{align}
        \hom{\bilinear{h}}(\Phi_H,\Psi_H) & = \sum_{K \in \mathcal{T}_H} \sum_{j= 1}^{J} \omega_K^j\, \homdisCoeKj\nabla \Phi_H \cdot  \nabla \Psi_H
        \label{hmm::bilinearform_dis}
    \end{align}
    with the discrete coefficient
    \begin{align}
        \label{hmm::operatordis}
        \homdisCoeKj\nabla \Phi_H \cdot  \nabla \Psi_H & =   \dfrac{1}{|	\cell{\delta}|}\int_{	\cellKj{\delta}} 	  \Coeffeps(x) \nabla \phihj(x)  \cdot \nabla  \psihj(x) \d{x}.
    \end{align}
    \begin{Bemerkung}
        In the discrete formulation \eqref{hmm::cellproblem_dis}, \eqref{hmm::operatordis} we replaced $a\Bigl(x^j_K, \dfrac{x}{\epsilon} \Bigr)$  by $\Coeffeps(x)$ since the structured form may be unknown in applications.
    \end{Bemerkung}
    \subsection{Error bounds for the fully discrete FE-HMM}
In this section we want to apply the abstract spatial error bounds to the discrete homogenized problem. 
\begin{Lemma}
    Let Assumption $\ref{model::assumption}$ hold true. Then the FE-HMM satisfies Assumptions $\ref{Spacedis::assumptions}$ and $\ref{Spacedis::Operators}$. 
\end{Lemma}
\begin{proof}
    $\,$ We set \[(\cdot,\cdot)_{\hildisnorm} = (\cdot,\cdot)_\hilnorm\quad, \quad  \energydisnorm =  \hom{\bilinear{h}} . \] 
    Hence, $\ritzH$ is the $L^2$-projection on $ \hildisSpac= (\disSpace, (\cdot,\cdot)_{\hildisnorm})$ with $\constdisH = 1$. Using \cite[Lemma 5]{Abdulle2009}, $\hom{\bilinear{h}}$ is coercive with constant $\lambda_H$ and bounded and therefore  Assumption \ref{Spacedis::assumptions}(i) is fulfilled, where $\constspdis =  0$ and $\constembdis = \lambda_H \constpoincare$. Then $\ritzV$ exists by the  Theorem of Lax-Migram and $\constdisV = \sqrt{\frac{\constbounded}{\lambda_H}}$. Consequently, Assumption \ref{Spacedis::Operators}(i)(iii) is fulfilled.  We set \[\bilinearbeta{H} = \interpol \bilinearbeta{} \quad,\quad \nonlindis = \interpol \nonlin \quad,\quad f_H = \interpol f,\]
   where we denote as $\interpol $ the nodal interpolation  operator as mapping from $H^2(\Omega) $ to $\soldisSpac =  ( \disSpace, \energydisnorm)$. Then  Assumption \ref{Spacedis::Operators}(ii) is fulfilled for a constant $\constinterpol$ due to standard interpolation error estimates for each subspace $\refSpac$ of $C^0(\overline{\Omega})$.
   Assumption \ref{Spacedis::assumptions}(ii),(iv)  follows directly from Assumption \ref{model::assumption} and boundedness of the interpolation operator.  In addition, for all $v_H,w_H \in \hildisSpac$ with $\|v_H\|_{\hildisnorm}, \|w_H\|_{\hildisnorm} \leq \rho_H$ we obtain 
   \begin{equation*}
       \label{Spaceerror::lipschitz}
       \begin{aligned}
           \| \nonlindis(\cdot,v_H)-\nonlindis(\cdot,w_H)\|_{\hildisnorm} \leq \constlipschitzdis \|v_H-w_H\|_{\hildisnorm}  
       \end{aligned}
   \end{equation*}
   and therefore \ref{Spacedis::assumptions}(iii) is fulfilled. 
\end{proof}
Next, we specify the required regularity of the exact solution in order to derive a priori error bounds.
\begin{Annahme}
    \label{Spaceerror::assumption}
    \phantom{-} \hspace*{0.1cm}
    \begin{enumerate}
        \item[(i) ] Let $C_1,C_2>0$ be constants independent of $\veps$ and $\delta$. For the small-scale parameter we assume \[ \varepsilon \leq C_1\delta^2 \leq C_2 H. \]
        \item[(ii) ] For the  discrete initial conditions we assume \[ \|\homdisnull{u} -\interpol \homnull{u}\|_{\energytdisnorm} + \| \homdisnull{v} - \interpol \homnull{v}\|_{\hildisnorm} \leq C H^p.\]
        \item[(iii) ] The right-hand side and the coefficients satisfy \[f \in \L{\infty}{\Hspace{\max\{2,p\}}}, \quad  \beta \in W^{p +1   , \infty}(\Omega), \quad a \in W^{1, \infty}(\Omega, W_{\#}^{1,\infty}(Y) ).\]
        \item[(iv) ]
        Let $T>0$ and $p \geq 1$. Then  we assume that the solution of the homogenized problem \eqref{model::HomLösunglinear1} satisfies the following regularity assumptions
        \begin{equation*}
            \begin{aligned}
                \hom{u} &\in \L{\infty}{\!\Hspace{p+1}}, \\  \partial_t \hom{u} &\in \L{\infty}{\!\Hspace{p+2}}, \\    \partial_{tt} \hom{u} &\in \L{\infty}{\!\Hspace{\max\{2,p\}}}.
            \end{aligned}
        \end{equation*}
        We  set
        \begin{equation*}
            \begin{aligned}
                \rho = 3 \max \Bigl\{1, \sqrt{\frac{\constbounded}{\lambda_H}}\Bigr\} \Bigl(\|\hom{u}\|_{\L{\infty}{\solSpac}} &+ \|\partial_t \hom{u}\|_{\L{\infty}{\hilSpac}}\Bigr)   + \constinterpol\|\partial_t \hom{u}\|_{\L{\infty}{\Hohnespace{\max\{2,p\}}}}
            \end{aligned}
        \end{equation*}
        Further, for the nonlinearity we assume
        \[
        \nonlin(x, \partial_t \hom{u}) \in  \L{\infty}{\Hspace{p}}.      \asseqqed
        \]
    \end{enumerate}
\end{Annahme}
We now estimate the errors terms occurring in Theorem \ref{fulldis::TheoremIMEX} in our specific setting.
\begin{Lemma}
    \label{Spaceerror::lemma:errorbound}
    Under the Assumption $\ref{model::assumption} $ and $\ref{Spaceerror::assumption}$, the spatial discretization error bounds of FE-HMM defined in Theorem $\ref{fulldis::TheoremIMEX}$ can be bounded by
    \[ \abstractspaceError{i} \leq C\bigl(H^p + \Bigl(\frac{h}{\veps}  \Bigr)^{2q} +  e_{\text{\emph{mod}}}\bigr).\]
    for $i = 1, \dots 4$, where $e_{\text{mod}}$ is a modeling error. The constant $C$ does not depend on $H$, $h$, and $\veps$.
\end{Lemma}
\begin{proof}
    $\,$ The terms $\abstractspaceError{i}$ for $ i = 1,3$ can be bounded using standard interpolation bounds (see ,e.g., \cite[Corollary 4.4.20]{Brenner2008Scott}) by \[\abstractspaceError{i} \leq C H^p. \] 
    For the errors of the nonlinearity $\nonlin$ we obtain by local Lipschitz continuity
    \begin{equation*}
        \begin{aligned}
            \|\ritzH \nonlin(\partial_t \hom{u}) - \nonlindis(\interpol \partial_t \hom{u}) \|_{\hildisnorm} \leq   C L_{\rho }H^p (  \| \nonlin( \partial_t \hom{u}) \|_{H^p} + \|\partial_t \hom{u}\|_{H^p}).
        \end{aligned}
    \end{equation*} 
    Further, for the strong damping term we obtain
    \begin{equation*}
        \begin{aligned}
            \|\ritzH \assdampOp \partial_t \hom{u} - \assdisdampOp \interpol \partial_t \hom{u} \|_{\hildisnorm} \leq C H^p(\|\beta \|_{W^{p+1,\infty}(\Omega)} + \|\beta \|_{C(\overline{\Omega})} + \|\partial_t \hom u\|_{H^{p+2}(\Omega)}) 
        \end{aligned}
    \end{equation*} 
    for $ w_H \in \disSpace$. 
    We obtain 
    \begin{equation*}
        \begin{aligned}
            \abstractspaceError{4} \leq C H^p. 
        \end{aligned}
    \end{equation*}
    It remains to bound the conformity error of the scalar products in $ \abstractspaceError{2}$. Let $\Phi_H , \Psi_H \in \soldisSpac$. Then 
    \begin{align*}
        |\bilinearhom{}	(\Phi_H,\Psi_H)- \hom{\bilinear{h}}(\Phi_H,\Psi_H)|  \leq \bigl( \sup_{\substack{ K \in \mathcal{T}_H \\ x_{K_j}\in K}} \| \homCoeKj(x_{K_j}) - \homdisCoeKj(x_{K_j})\|_F + {e_{\text{qf}}}\bigr) \|\Phi_H\|_{\solSpac} \|\Psi_H\|_{\solSpac}, 
    \end{align*}
    where $\|\cdot\|_F$ is the Frobenius norm and $e_{\text{qf}}$  is a quadrature error arising from the quadrature in \eqref{hmm::bilinearform}. 
    For estimating this we define the auxiliary problem:
    Find $\phihbarj - \Phi^{lin}_{H,K,j}  \in  W(\cellKj{\delta})$ such that
    \begin{align}
        \int_{\cellKj{\delta}}\Coeffeps (x) \nabla \phihbarj \cdot \nabla z \d{x}= 0
    \end{align}
    for all $z \in W(\cellKj{\delta})$. 
    Define
    \begin{align*}
        \hilfAj\nabla \Phi_H \cdot \nabla \Psi_H & =   \dfrac{1}{|\cell{\delta}|}\int_{\cellKj{\delta}} 	  \Coeffeps(x) \nabla  \phihbarj \cdot \nabla \Psi_H \d{x}. \\
    \end{align*}
    We now split the micro error into two terms
    \begin{align*}
        \| \homCoeKj(x_{K_j}) - \homdisCoeKj(x_{K_j})\|_F \leq \| \homCoeKj(x_{K_j}) - \hilfAj(x_{K_j})\|_F + \| \hilfAj(x_{K_j}) - \homdisCoeKj(x_{K_j}) \|_F.
    \end{align*}
    The first term represents the model error $e_{\text{mod}}$ caused by the use of wrong cell sizes or wrong boundary conditions.
    The second term represents the discretization error of the cell problems and can be estimated as follows using \cite[Lemma 4.7]{Abdulle2012Engquist} and Remark \ref{remark::qf}:
    \begin{align*}
        \| \hilfAj(x_{K_j}) - \homdisCoeKj(x_{K_j}) \|_F + {e_{\text{qf}}} \leq C \Bigl(\frac{h}{\epsilon}\Bigr)^{2q}.
    \end{align*}
\end{proof}
\begin{Bemerkung}
    The model error $e_{\text{mod}}$ can be bounded depending on the type of boundary conditions selected. 
    If $\Coeffeps$ is sufficiently regular 
    \begin{align*}
        e_{\text{mod}} \leq  C
        \begin{cases}
             \delta, \qquad &\text{if } \frac{\delta}{\epsilon} \in \N  \text{ and  periodic b.c.,}\\
             \delta + \frac{\veps}{\delta}, \qquad &\text{else. } 
            \end{cases} 
    \end{align*}
    The proofs can be found in \cite[Proposition 14]{Abdulle2009} and \cite[Theorem 3.4]{Du2010Ming}. Note that  \[ W^{1,\infty}(\Omega , W^{1,\infty}_{\#}(Y) ) \hookrightarrow C^{0,1}(\Omega , W^{1,\infty}_{\#}(Y))\] (see \cite[Theorem 5.2]{Arendt2018Kreuter}). 
\end{Bemerkung}
We now apply the abstract time error bounds to the homogenized equation \eqref{model::HomLösunglinear1}.
\begin{Satz}
    \label{fullerror::IMEX}
    Let Assumption $\ref{model::assumption}$ be fulfilled and assume that the solution $\hom{u}$ of \eqref{model::HomLösunglinear1} satisfies
    \[u \in C^4([0,T]; \hilSpac) \cap C^3([0,T];\solSpac)\cap C^2([0,T];H^2(\Omega)). \]
    \begin{itemize}
        \item[(a) ](Semi-discrete): There exists $\tau^\star>0$, such that for all $\tau < \tau^\star$ the error of the \timeRscheme applied on \eqref{model::HomLösunglinear1} can be bounded by 
        \begin{equation*}
            \begin{aligned}
                \| \homn{u} - \hom{u}(t_n)\|_{\solSpac} + \| \homn{v} - \partial _t\hom{u}(t_n) \|_{\hilSpac} \leq C t_n \e^{M  t_n  }  \tau^2,
            \end{aligned}
        \end{equation*} 
        where $M \coloneqq \bigl(\frac{1}{2}{\constcoercive} \constpoincare + \constlipschitz (2+ \tauhalb \constlipschitz)\bigr) $.
        The constant $C$ depends only on the derivatives of $\hom{u}$ and  $ \constlipschitz$, $T$.
        \item[(b) ](Full-discrete): Let Assumption $\ref{Spaceerror::assumption}$ be fulfilled. Then there exist $ \tau^\star, h^\star>0$ such that for all $ \tau <\tau^\star$ and  $h < h^\star $  the errors of the fully discrete solution of the \timeRscheme and the FE-HMM satisfy the following bounds
        \begin{align*}
            \| \homdisn{u} - \hom{u}(t_n)\|_{\solSpac} \!+ \!\| \homdisn{v} - \partial_t\hom{u}(t_n) \|_{\hilSpac} \leq C t_n \e^{M t_n  }  \bigl(\tau^2  \!+ \! H^p +  \Bigl(\frac{h}{\veps}  \Bigr)^{2q}  \!\!+ \!  e_{\text{\emph{mod}}}\bigr),
        \end{align*}
        where $M \coloneqq \bigl(\frac{1}{2}{\lambda_H} \constpoincare + \constlipschitz (2+ \tauhalb \constlipschitz)\bigr)$.
        The constant $C$ depends only on the derivatives of $\hom{u}$ and  $ \constlipschitz$, $T$.
    \end{itemize} 
\end{Satz}
\section{Numerical Examples \label{num_example}}
 In the following we present numerical results to demonstrate the accuracy of the error estimates. For this, we use an implementation in the  $C\!++$ FEM-library \textit{deal.II}. We verify numerically the space and time convergence rates of Theorem \ref{fullerror::IMEX} for periodic coefficients. 
In the following  we define the error function for all $ \tau > 0 $ as \[ \mathbf{E}_\tau(t_n) \coloneqq \dfrac{\|\nohomdisn{u}
 - \nohom{u}(t_n)\|_{H^1(\Omega)} + \|\nohomdisn{v}- \nohom{u}^\prime(t_n) \|_{L^2(\Omega)} }{\|
  \nohom{u}(t_n)\|_{H^1(\Omega)} + \| \nohom{u}^\prime(t_n) \|_{L^2(\Omega)}}, \quad \text{where }  t_n = n\tau. \] Note that in our example the denominator will never be zero. 
\subsection{Implementation}
To explain the implementation we first introduce some notation. We denote by  $\solvecfem,\solprimevecfem \in \R^{N_H}$ the coefficient vector corresponding to the solution $\nohomdis{u},\nohomdis{v} \in \soldisSpac$. Furthermore, we set $\massmatrix \in \R^{N_H \times N_H}$ as the mass matrix. We denote $\stiffmatrix, \dampmatrix \in \R^{N_H\times N_H}$ as the stiffness matrices and $\nonlinearvec, \loadvec \in \R^{N_H}$ as the load vectors to the discrete quantities $\assdisOp, \assdisdampOp$ and $\nonlindis, f_H$. In the following, we set \[\nonlinearvec^n(\nohomn{\solvecfem},\nohomn{\solprimevecfem})  = \nonlinearvec(\nohomn{\solvecfem},\nohomn{\solprimevecfem}) + \loadvec^{n}. \] From \eqref{errorbounds::IMEXR1} we obtain for the half step the following second-order formation
\begin{subequations}
    \label{Implementation::generalhalf}
    \begin{alignat}{2}
        \massmatrix\nohomhalbe{\solprimevecfem} &= \massmatrix\nohomn{\solprimevecfem}  - \frac{\tau}{2} \stiffmatrix \nohomn{\solvecfem} - \frac{\tau^2}{4} \stiffmatrix  \nohomhalbe{\solprimevecfem} - \frac{\tau}{2} \dampmatrix\nohomhalbe{\solprimevecfem} + \frac{\tau}{2}  \nonlinearvec^n(\nohomn{\solvecfem},\nohomn{\solprimevecfem}), \label{Implementation::generalhalf2} \\
        \nohomhalbe{\solvecfem} &= \nohomn{\solvecfem} + \tauhalb \nohomhalbe{\solprimevecfem}. \label{Implementation::generalhalf1}   
    \end{alignat}
\end{subequations}
Using this, we can calculate $\nohomnp{\solprimevecfem}$ and $\nohomnp{\solvecfem}$ for the full-step \eqref{errorbounds::IMEXR2} by 
\begin{subequations}
    \label{Implementation::generalfull}
    \begin{alignat}{2}
        \massmatrix\nohomnp{\solprimevecfem} &= 2\massmatrix\nohomhalbe{\solprimevecfem} - \massmatrix\nohomn{\solprimevecfem} + \tau \bigl( \nonlinearvec^{n + \frac{1}{2}}(\nohomhalbe{\solvecfem},\nohomhalbe{\solprimevecfem})   - \nonlinearvec^n(\nohomn{\solvecfem},\nohomn{\solprimevecfem}) \bigr), \label{Implementation::generalfull1} \\
        \nohomnp{\solvecfem} &= \nohomn{\solvecfem} + \tau\nohomhalbe{\solprimevecfem}, \label{Implementation::generalfull2}
    \end{alignat}
\end{subequations}
In total, we have to solve two linear systems \eqref{Implementation::generalhalf2}, \eqref{Implementation::generalfull1} in each time step with one application of $\stiffmatrix$, one of $\dampmatrix$ and two evaluations of the nonlinearity $\nonlinearvec$.  \\
In the special case $\nonlinearvec^n(\nohomn{\solvecfem},\nohomn{\solprimevecfem}) = \nonlinearvec^n(\nohomn{\solvecfem})$ only one linear system, namely \eqref{Implementation::generalhalf2}, must be solved, since we do not need to calculate $\nohomnp{\solprimevecfem}$, as only $\massmatrix\nohomnp{\solprimevecfem}$ is required.
\subsection{Setting}
We set as our domain $\Omega = (0,1)^2$, $T = 1$ and consider the homogenized nonlinear wave equation 
\begin{equation}
    \begin{aligned}
        \label{example::equation}
        \partial_{tt} \hom{u} &- \nabla\cdot (\hom{a}(x)\nabla \hom{u}) - 0.01\Delta (\partial_t \hom{u}) \\ &+ \sgn(\partial_{t}\hom{u})(|\partial_{t}\hom{u} + 0.0001|^{0.6} - 0.0001^{0.6}) = f 
    \end{aligned}
\end{equation}
with homogeneous Dirichlet boundary conditions. 
For the microscopic problem we consider periodic boundary conditions with $\veps = 2 ^{-15}$ and  $\delta= 2^{-13}$, where we use $\mathcal{P}_1$ linear elements on a uniform microscopic mesh with grid width $h$.  We set
\begin{align*}
    a^{\epsilon}(x_1,x_2) =  \begin{pmatrix}
        0.33   + 0.15(\sin(2\pi x_1) + \sin(2\pi \frac{x_1}{\epsilon})) & 0                                                               \\
        0                                                                & 0.33   + 0.15(\sin(2\pi x_1) + \sin(2\pi \frac{x_1}{\epsilon})) \\
    \end{pmatrix}
\end{align*}
as coefficient. According to \cite[Chapter 5.3]{Abdulle2014}, the homogenized coefficient can be calculated exactly as
\begin{equation}
    \label{example::exact_coefficient}
    \begin{aligned}
        a^{hom}(x_1,x_2) =
         \begin{pmatrix}
            0.3\cdot\sqrt{(1.1+0.5\cdot \sin(2\pi x_1))^2-0.25} & 0                                  \\
            0                                                 & 0.3\cdot(1.1+0.5\cdot\sin(2\pi x_1))
        \end{pmatrix} .
    \end{aligned}
\end{equation}
As an exact solution of the homogenized equation we choose
\begin{align*}
    \hom{u}(t,x_1,x_2) = \e^{\pi t} \sin(\pi x_1^2) \sin(\pi x_2^2)
\end{align*}
with a corresponding right-hand side and initial conditions. For the macroscopic discretization we use finite elements on a uniform quadrilateral mesh with meshsizes $ H = 2^{-k}$. 
\subsection{Results}
In Figure \ref{example::space} we consider the spatial error for different microscopic meshsizes for a fixed time step-size $\tau = 10^{-3}$. Here we use finite elements of order $p=1$ as macroscopic space discretization. As predicted in Theorem \ref{fullerror::IMEX}, we obtain first order convergence till we reach the plateau of the space discretization of the cell problems. 
\begin{figure}
    \centering
    \subcaptionbox{\label{example::space}}{%
    \vspace*{-0.095cm}
      \input{plots/errorplots/Example0/space/TikzIMEXRerror_total_e15d13-damped.tikz}\hspace*{1em}%
    } \hspace*{-0.3cm}
    \subcaptionbox{\label{example::time}}{%
      \input{plots/errorplots/Example0/time/Tikzerror_total-exact-damped.tikz}\hspace*{1em}%
    }
    \caption{The error $\errorfunc(1)$ of finite element heterogeneous multiscale method plotted against the macroscopic mesh size $H$ on the left and the time step-size $\tau$ on the right. \label{example::figure} }
\end{figure}
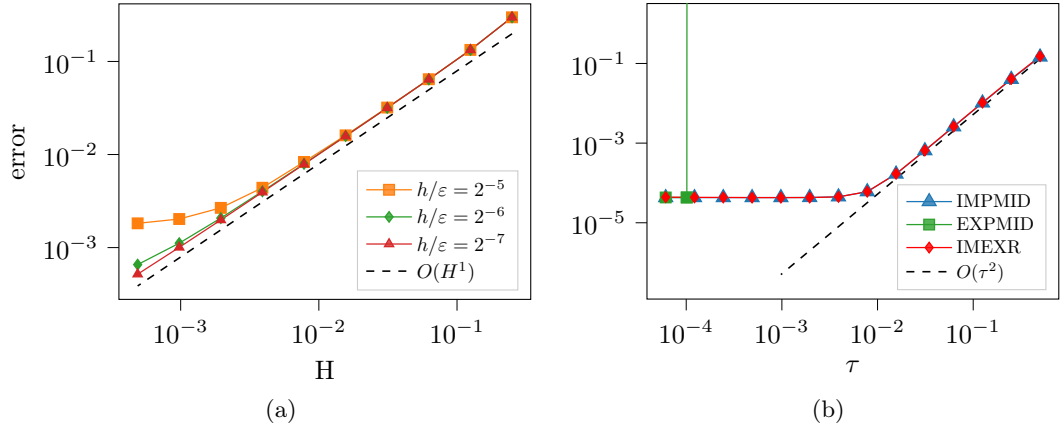
In Figure \ref{example::time} we compare the implicit and explicit midpoint rule  with the IMEX scheme. 
For the macroscopic space discretization, we use piecewise quadratic elements with a fixed $H = 2^{-7}$.  In addition, we use the formula \eqref{example::exact_coefficient} for the homogenized coefficient  here, since we are focusing on the time discretization error. 
One can observe the predicted convergence rate for the IMEX scheme as well as the second-order convergence for the implicit midpoint rule until the spatial discretization error is reached. The explicit midpoint rule is only stable under a strong CFL condition and the error immediately reaches the plateau of the spatial discretization error.
\paragraph*{Conflict of interest} The authors declare that they have no conflict of interest to this work.
\paragraph*{Data availability} The codes to reproduce the experiments are available  on \url{https://doi.org/10.35097/gkHSEiWVXPBkLPxL}.
\bibliographystyle{abbrv}
\appendix%
\bibliography{source}
	
\end{document}

%% file: plots/errorplots/Example0/space/TikzIMEXRerror_total_e15d13-damped.tikz
\begin{tikzpicture}

\definecolor{crimson2143940}{RGB}{214,39,40}
\definecolor{darkgray176}{RGB}{176,176,176}
\definecolor{darkorange25512714}{RGB}{255,127,14}
\definecolor{forestgreen4416044}{RGB}{44,160,44}
\definecolor{lightgray204}{RGB}{204,204,204}
\definecolor{steelblue31119180}{RGB}{31,119,180}

\begin{axis}[
height=5.5cm,
legend cell align={left},
legend style={
  nodes={scale=0.75, transform shape},
  fill opacity=0.8,
  draw opacity=1,
  text opacity=1,
  at={(0.97,0.03)},
  anchor=south east,
  draw=lightgray204
},
log basis x={10},
log basis y={10},
tick align=outside,
tick pos=left,
width=7.0cm,
x grid style={darkgray176},
xlabel={H},
xmin=0.000357442796861725, xmax=0.341510064188599,
xmode=log,
xtick style={color=black},
xtick={1e-05,0.0001,0.001,0.01,0.1,1,10},
xticklabels={
  \(\displaystyle {10^{-5}}\),
  \(\displaystyle {10^{-4}}\),
  \(\displaystyle {10^{-3}}\),
  \(\displaystyle {10^{-2}}\),
  \(\displaystyle {10^{-1}}\),
  \(\displaystyle {10^{0}}\),
  \(\displaystyle {10^{1}}\)
},
y grid style={darkgray176},
ylabel={error},
ymin=0.000277651479289042, ymax=0.415758691149736,
ymode=log,
ytick style={color=black},
ytick={1e-05,0.0001,0.001,0.01,0.1,1,10},
yticklabels={
  \(\displaystyle {10^{-5}}\),
  \(\displaystyle {10^{-4}}\),
  \(\displaystyle {10^{-3}}\),
  \(\displaystyle {10^{-2}}\),
  \(\displaystyle {10^{-1}}\),
  \(\displaystyle {10^{0}}\),
  \(\displaystyle {10^{1}}\)
}
]

\addplot [line width=0.48pt, darkorange25512714, mark=square*, mark size=2, mark options={solid}]
table {%
0.25 0.2982
0.125 0.1331
0.0625 0.06442
0.03125 0.03197
0.015625 0.01608
0.0078125 0.008327
0.00390625 0.004395
0.001953125 0.002673
0.0009765625 0.002028
0.00048828125 0.001823
};
\addlegendentry{$h/\varepsilon = 2^{-5}$}
\addplot [line width=0.48pt, forestgreen4416044, mark=diamond*, mark size=2, mark options={solid}]
table {%
0.25 0.2982
0.125 0.133
0.0625 0.06428
0.03125 0.03172
0.015625 0.01581
0.0078125 0.007938
0.00390625 0.004025
0.001953125 0.002075
0.0009765625 0.001116
0.00048828125 0.0006564
};
\addlegendentry{$h/\varepsilon = 2^{-6}$}
\addplot [line width=0.48pt, crimson2143940, mark=triangle*, mark size=2, mark options={solid}]
table {%
0.25 0.2982
0.125 0.133
0.0625 0.06425
0.03125 0.03168
0.015625 0.01575
0.0078125 0.007861
0.00390625 0.003939
0.001953125 0.001981
0.0009765625 0.001006
0.00048828125 0.0005181
};
\addlegendentry{$h/\varepsilon = 2^{-7}$}
\addplot [semithick, black, dashed]
table {%
0.25 0.1982
0.00048828125 0.000387109375
};
\addlegendentry{$O(H^{1})$}
\end{axis}

\end{tikzpicture}

%% file: plots/errorplots/Example0/time/Tikzerror_total-exact-damped.tikz
\begin{tikzpicture}

\definecolor{darkgray176}{RGB}{176,176,176}
\definecolor{darkorange25512714}{RGB}{255,127,14}
\definecolor{forestgreen4416044}{RGB}{44,160,44}
\definecolor{lightgray204}{RGB}{204,204,204}
\definecolor{steelblue31119180}{RGB}{31,119,180}

\begin{axis}[
height=5.5cm,
legend cell align={left},
legend style={
  nodes={scale=0.7, transform shape},
  fill opacity=0.8,
  draw opacity=1,
  text opacity=1,
  at={(0.97,0.03)},
  anchor=south east,
  draw=lightgray204
},
log basis x={10},
log basis y={10},
tick align=outside,
tick pos=left,
width=7.0cm,
x grid style={darkgray176},
xlabel={\(\displaystyle \tau\)},
xmin=3.88963989654698e-05, xmax=0.78458419832365,
xmode=log,
xtick style={color=black},
xtick={1e-06,1e-05,0.0001,0.001,0.01,0.1,1,10},
xticklabels={
  \(\displaystyle {10^{-6}}\),
  \(\displaystyle {10^{-5}}\),
  \(\displaystyle {10^{-4}}\),
  \(\displaystyle {10^{-3}}\),
  \(\displaystyle {10^{-2}}\),
  \(\displaystyle {10^{-1}}\),
  \(\displaystyle {10^{0}}\),
  \(\displaystyle {10^{1}}\)
},
y grid style={darkgray176},
ymin=1.21402058109505e-07, ymax=3.19959454,
ymode=log,
ytick style={color=black},
ytick={1e-09,1e-07,1e-05,0.001,0.1,10,1000,100000,10000000,1000000000},
yticklabels={
  \(\displaystyle {10^{-9}}\),
  \(\displaystyle {10^{-7}}\),
  \(\displaystyle {10^{-5}}\),
  \(\displaystyle {10^{-3}}\),
  \(\displaystyle {10^{-1}}\),
  \(\displaystyle {10^{1}}\),
  \(\displaystyle {10^{3}}\),
  \(\displaystyle {10^{5}}\),
  \(\displaystyle {10^{7}}\),
  \(\displaystyle {10^{9}}\)
}
]
\addplot [line width=0.48pt, steelblue31119180, mark=triangle*, mark size=3, mark options={solid}]
table {%
0.5 0.1453
0.25 0.03991
0.125 0.01022
0.0625 0.00257
0.03125 0.0006451
0.015625 0.0001672
0.0078125 6.024e-05
0.00390625 4.48e-05
0.001953125 4.31e-05
0.000976562 4.288e-05
0.000488281 4.286e-05
0.000244141 4.286e-05
0.00012207 4.286e-05
6.1035e-05 4.286e-05
};
\addlegendentry{IMPMID}
\addplot [line width=0.48pt, forestgreen4416044, mark=square*, mark size=2, mark options={solid}]
table {%
0.000102041 959100
0.000101523 4.318e-05
0.00010101 4.318e-05
6.1035e-05 4.317e-05
};
\addlegendentry{EXPMID}
\addplot [line width=0.48pt, red, mark=diamond*, mark size=2, mark options={solid}]
table {%
0.5 0.1519
0.25 0.04095
0.125 0.01043
0.0625 0.002618
0.03125 0.0006566
0.015625 0.0001698
0.0078125 6.04e-05
0.00390625 4.468e-05
0.001953125 4.305e-05
0.000976562 4.287e-05
0.000488281 4.293e-05
0.000244141 4.322e-05
0.00012207 4.341e-05
6.1035e-05 4.399e-05
};
\addlegendentry{IMEXR}
\addplot [semithick, black, dashed]
table {%
0.5 0.1309
0.000976562 4.99343360742318e-07
};
\addlegendentry{$O(\tau^{2})$}
\end{axis}

\end{tikzpicture}